\documentclass{amsart}
\usepackage{ifthen}
\usepackage{amsmath,amssymb,latexsym}
\usepackage{times}
\usepackage[T1]{fontenc}

\font\goth=eusm10
\newcommand\F{\mathcal F}
\newcommand\E{\mathcal E}
\newcommand\N{\mathcal N}
\newcommand\T{\mathcal T}
\newcommand\La{\mathcal L}
\newcommand\G{\mathcal G}
\newcommand\V{\mathcal V}

\newcommand\Pmc{\mathcal P}
\newcommand\Ii{\hbox{\goth I}}
\newcommand\Oc{\hbox{\goth O}}

\newcommand\FF{\mathbb{F}}
\newcommand\ZZ{\mathbb{Z}}
\newcommand\CC{\mathbb{C}}
\newcommand\NN{\mathbb{N}}
\newcommand\RR{\mathbb{R}}
\newcommand\QQ{\mathbb{Q}}

\newcommand\Pp{\mathbb P}
\newcommand\PR{\mathbb{P}^r}
\newcommand\Pd{\mathbb{P}^2}
\newcommand\Pt{\mathbb{P}^3}

\numberwithin{equation}{section}

\newtheorem{theorem}[equation]{Theorem}
\newtheorem{coroll}[equation]{Corollary}

\newtheorem{lemma}[equation]{Lemma}
\newtheorem{proposition}[equation]{Proposition}

\theoremstyle{definition}
\newtheorem{definition}[equation]{Definition}
\newtheorem{notation}[equation]{Notation}

\newtheorem{remark}[equation]{Remark}

\begin{document} 

\title{Families of nodal curves on projective 
threefolds and their regularity via postulation of nodes}
\author{Flaminio Flamini}

\email{flamini@matrm3.mat.uniroma3.it}
\curraddr{Dipartimento di Matematica, Università degli Studi de L'Aquila, Via 
Vetoio - Coppito 1, 67100 L'Aquila - Italy}

\thanks{The author is a member of Cofin GVA, EAGER and GNSAGA-INdAM}

\begin{abstract}{The main purpose of this paper is to introduce a 
new approach to study families of nodal curves on projective threefolds. 
Precisely, given $X$ a smooth projective threefold, 
$\E$ a rank-two vector bundle on $X$, $L$ a very ample line bundle on $X$ and 
$k \geq 0$, $\delta >0 $ integers and 
denoted by ${\V}_{\delta} ({\E} \otimes L^{\otimes k})$ the subscheme of 
${\Pp}(H^0({\E} \otimes L^{\otimes k}))$ 
parametrizing global sections of ${\E} \otimes L^{\otimes k}$ whose zero-loci are 
irreducible and $\delta$-nodal curves on $X$, we present a new cohomological 
description of the tangent 
space $T_{[s]}({\V}_{\delta} ({\E} \otimes L^{\otimes k}))$ at a point
$[s]\in {\V}_{\delta} ({\E} \otimes L^{\otimes k})$. This description enable us to
determine effective and uniform upper-bounds 
for $\delta$, which are linear polynomials in $k$, such that the family 
${\V}_{\delta} ({\E} \otimes L^{\otimes k})$ is smooth 
and of the expected dimension ({\em regular}, for short). 
The almost-sharpness of our bounds is shown by some interesting examples. 
Furthermore, when $X$ is assumed to be a Fano or a Calaby-Yau 
threefold, we study in detail 
the regularity property of a point $[s] \in {\V}_{\delta} ({\E} \otimes L^{\otimes k})$
related to the postulation of the nodes of its zero-locus $C = V(s) \subset X$.
Roughly speaking, when the nodes of $C$ are assumed to be in 
general position either on $X$ or on an irreducible divisor of $X$ 
having at worst log-terminal singularities or to lie on a l.c.i. and subcanonical curve 
in $X$, 
we find upper-bounds on $\delta$ which are, respectively, cubic, quadratic and linear 
polynomials in $k$ ensuring the regularity of ${\V}_{\delta} ({\E} \otimes L^{\otimes k})$ at 
$[s]$. Finally, when $X= \Pt$, we also discuss some interesting geometric properties of the
curves given by sections parametrized by ${\V}_{\delta} ({\E} \otimes \Oc_X(k))$.}
\end{abstract}

\maketitle

\section*{Introduction}\label{S:intro}
The theory of families of singular curves with fixed 
invariants (e.g. geometric genus, singularity type, 
number of irreducible components, etc.) and which are contained in a 
projective variety $X$ has been extensively 
studied from the beginning of Algebraic Geometry and it 
actually receives a lot of attention, partially due to 
its connections with other fields of geometry and physics. 
Indeed, the interest about this arguments has grown essentially for two reasons: on the 
one hand, the theory of strings of nuclear physicists deals with enumerative geometry 
for rational curves contained in some projective threefolds; on the other hand, 
the study of singular curves is naturally related with the hyperbolic geometry of 
complex projective varieties.

Nodal curves play a central role in the subject of singular curves. 
Families of irreducible and $\delta$-nodal curves on a given projective variety 
$X$ are usually 
called {\em Severi varieties} of irreducible, $\delta$-nodal curves in $X$. 
The terminology "Severi variety" is due to the classical case of families of 
nodal curves on $X= \Pd$, which was first studied by Severi (see \cite{Sev}).

The case in which $X$ is a smooth projective surface has recently given rise 
to a huge amount of literature (see, for example, \cite{CH}, \cite{CC}, \cite{CL}, 
\cite{CS}, \cite{F1}, \cite{GLS}, \cite{Harris}, \cite{Ran}, \cite{S}, \cite{Xu} just to 
mention a few. For a detailed chronological 
overview, the reader is referred for example to Section 2.3 in \cite{F2} and 
to its bibliography). This depends not only on the great interest in the subject, but 
also because for a Severi 
variety $V$ on an arbitrary projective variety $X$ there are several problems concerning $V$ 
like non-emptyness, smoothness, irreducibility, dimensional computation as well as 
enumerative and moduli properties of the family of curves 
it parametrizes.

On the contrary, in higher dimension only few results are known. Precisely, 
on the one hand we have non-emptyness and 
enumerative results for some classes 
of varieties, which are relevant for applications; on the other hand, 
some other results of non-emptyness and smoothness are given 
only for families of nodal curves in 
projective spaces (see e.g. \cite{BC} and \cite{S}).

Therefore, we feel some lack of systematic 
studies for what should be the next relevant case, from the point of view 
of Algebraic Geometry: families of nodal curves on projective threefolds.

The purpose of this paper is twofold: first, we introduce a new 
method to determine when a given (non-empty) Severi variety 
on a smooth projective threefold $X$ is smooth and of the expected dimension 
- {\em regular}, for short - at a given point (for details, see 
Definition \ref{def:0}); then, we apply this method to find geometric and 
numerical sufficient conditions for the regularity property of Severi varieties.

In general, a natural approach to the regularity problem is to use deformation theory of 
nodal curves in a smooth ambient variety. In the surface case, 
if $V_{\mid {\Oc}_S(D) \mid, \delta}$ denotes 
the Severi variety of irreducible and $\delta$-nodal curves 
in the linear system $|{\Oc}_S(D)|$ on a smooth projective surface $S$, 
it is well known that, when $V_{\mid {\Oc}_S(D) \mid, \delta} \neq \emptyset$, 
its expected codimension in $|{\Oc}_S(D)|$ is $\delta$; moreover, if 
$[C] \in V_{\mid {\Oc}_S(D)\mid, \; \delta}$ parametrizes 
a curve $C$ whose set of nodes is denoted by $\Sigma$, the 
Zariski tangent space at $[C]$ is  
$$ T_{[C]}(V_{\mid {\Oc}_S(D)\mid, \; \delta}) \cong 
H^0(S, \Ii_{\Sigma/S}\otimes {\Oc}_S(D))/< C >,$$where 
$\Ii_{\Sigma/S}$ denotes the ideal sheaf of $\Sigma$ in $S$ (see, 
for example, \cite{CS}). Thus, since the relative obstruction space is contained
in $H^1(S, \Ii_{\Sigma/S} \otimes {\Oc}_S(D))$, the regularity of 
$V_{\mid {\Oc}_S(D)\mid, \; \delta}$ 
at $[C]$ holds iff $\Sigma$ imposes independent conditions to $|{\Oc}_S(D)|$;
in particular, a sufficient 
condition for the regularity at $[C]$ is $h^1(S, \Ii_{\Sigma/S} \otimes 
{\Oc}_S(D))=0$.   
 
In the threefold case, one obtains a partially similar organization for 
curves which are zero-loci of sections of a rank-two vector bundle $\F$ 
on $X$, so that ${\Pp}(H^0(X, {\F}))$ (which plays the same role of 
$|{\Oc}_S(D)|$) somehow gives a projective space 
dominating a subvariety in which the curves move. As in the surface case,
if ${\V}_{\delta} ({\F})$ denotes the subscheme of ${\Pp}(H^0(X, {\F}))$ parametrizing 
global sections whose zero-loci are irreducible and $\delta$-nodal curves 
in $X$, with a little abuse of terminology we shall always use the term 
{\em Severi variety} to refer to ${\V}_{\delta} ({\F})$. In several cases 
- e.g. when $\F$ is a stable and aCM rank-two vector bundle on $X$ 
(for definitions, see \cite{CM}) - 
this is not an abuse, since ${\V}_{\delta} ({\F})$ actually parametrizes irreducible nodal 
curves on $X$ (see Lemma 4.3 in \cite{CM}).

\begin{notation}\label{not:1}
In the sequel, we write
$[s] \in {\V}_{\delta} ({\F})$ to intend that the global section $s \in H^0(X, \F)$
determines the corresponding point $[s]$ of the scheme ${\V}_{\delta} ({\F})$.
We also denote by $C_s$ (or simply $C$, when this
does not create ambiguity) the zero-locus of the given section $s$, i.e.
$C= V(s) \subset X$. 
\end{notation}

As in the surface case, when ${\V}_{\delta} ({\F})$ is not empty then
its expected codimension is $\delta$ (see Proposition \ref{prop:0}); however, if
$[s] \in {\V}_{\delta} ({\F})$ and if $\Sigma$ denotes the set of nodes of the corresponding
curve $C_s \subset X$, now we have
$$T_{[s]}({\V}_{\delta} ({\F})) \supset H^0(X, \Ii_{\Sigma/X}\otimes {\F})/< s >,$$so the latter
is the tangent space at $[s]$ to a subscheme of ${\Pp}(H^0(X, {\F}))$ of a higher
expected codimension.

We thus present a systematic study of equisingular deformation
theory for the elements parametrized by ${\V}_{\delta} ({\F})$ on $X$.
Precisely in the following result, which is the core of the paper,
we introduce a new cohomological description
of the tangent space $T_{[s]}({\V}_{\delta} ({\F}))$.

\vskip 8pt

\noindent
{\bf Theorem} (see Theorem \ref{prop:3.fundamental})
{\em Let $X $ be a smooth projective threefold. Let $\F$ be a globally generated rank-two vector
bundle on $X$ and let $\delta$ be a positive integer. Let
\[\begin{aligned}
{\V}_{\delta}({\F}) := & \{[s] \in {\Pp}(H^0({\F})) \; | \; C_s := V(s)
\subset X \; {\rm is \; irreducible \;} \\
 & {\rm  with \; only} \;  \delta \; {\rm nodes \; as \; singularities} \}.
\end{aligned}\]Consider $[s] \in {\V}_{\delta}({\F})$ and
denote by $\Sigma$ the set of nodes of the corresponding curve
$C_s \subset X$. Let$${\mathcal P} := {\Pp}_{X}({\F}) \stackrel{\pi}{\longrightarrow} X$$be
the
projective space bundle together with its natural projection $\pi$ on $X$ and denote
by ${\Oc}_{\mathcal P}(1)$ its tautological line bundle. Let $T^1_{C_s}$ be the 
{\em first cotangent sheaf} of $C_s$ (see \eqref{eq:T1}) and let
$$\Sigma^1 : = {\Pp}_X(T^1_{C_s}) \subset {\Pmc}$$be the
zero-dimensional subscheme of $\Pmc$ of length $\delta$, determined by the
surjection ${\F} \to   {T^1}_{C_s} \to 0$. Denote by $\Ii_{\Sigma^1/{\Pmc}}$ the ideal 
sheaf of $\Sigma^1$ in $\Pmc$. Then
\begin{itemize}
\item[(i)] $\Sigma^1$ is a set of $\delta$ rational double points for the
divisor $D_s \in |{\Oc}_{\mathcal P}(1)|$, corresponding to the given section
$s \in H^0(X, \F)$, and
\item[(ii)] the subsheaf of $\F$, defined by
$${\F}^{\Sigma} := \pi_* ({\Ii}_{\Sigma^1/{\Pmc}} \otimes {\Oc}_{\mathcal P}(1)),$$is such
that its global sections (modulo the
one dimensional subspace $<s>$) parametrize first-order deformations of
$s \in H^0(X, \F)$ which are equisingular.

\noindent
In particular, we have
$$\frac{H^0(X, {\F}^{\Sigma})}{< s >} \cong T_{[s]} ({\V}_{\delta}({\F})) \subset T_{[s]}
({\Pp}(H^0({\F}))) \cong \frac{H^0(X, {\F})}{< s >}.$$
\end{itemize}
}

\vskip 8pt

\noindent
We want to briefly remark that the above result can be the starting point for the characterization
of tangent spaces to such families on a smooth projective $n$-fold $Y$, with
$n \geq 4$. The main difference from the threefold case
is that one should work inside a suitable incidence variety
${\mathcal I} \subset {\Pp}_{Y}({\F}) \times {\Pp}_{Y}({\F}^{\vee})$, 
where $rank(\F) = n-1$.

By using the above characterization of $T_{[s]} ({\V}_{\delta}({\F}))$, we are able to determine geometric and 
numerical sufficient conditions for the regularity of ${\V}_{\delta}({\F})$ at the point 
$[s]$.

Indeed, we first prove the following:

\vskip 8pt

\noindent
{\bf Theorem} (see Theorem \ref{thm:9bis})
{\em Let $X $ be a smooth projective threefold. 
Let $\E$ be a globally generated rank-two vector bundle, $L$ be a very ample line bundle and 
$k \geq 0$ and $\delta >0 $ be integers. If $$(*) \;\;\; \delta \leq k+1,$$then 
${\V}_{\delta}({\E} \otimes L^{\otimes k})$ is smooth and of the expected dimension (i.e. regular) 
at each point.}

\vskip 8pt

\noindent
Therefore, the above result determine sufficient conditions in order that 
${\V}_{\delta}({\E} \otimes L^{\otimes k})$ is 
regular everywhere. Observe also that the bound $(*)$ is uniform, i.e. it does not depend 
on the postulation of the nodes of the curves related to the elements parametrized 
by ${\V}_{\delta}({\E} \otimes L^{\otimes k})$; furthermore, the bound is almost-sharp, 
as one can deduce from 
Example 3.2 in \cite{BC} and from our Remarks \ref{rem:17} and \ref{rem:36}. 
We also stress that the above result generalizes 
what proved in \cite{BC}, mainly 
because our approach more generally holds for families of nodal curves on smooth projective 
threefolds but also because, even in the case of $X = \Pt$, main subject of \cite{BC}, 
our bounds are effective and not asymptotic as Proposition 3.1 in \cite{BC}. 

After this, in \S 5 we focus on the case of $X \subset \PR$ either a Fano or a Calaby-Yau 
threefold, with $L = \Oc_X(1)$, and 
we "stratify" the regularity property in terms of the {\em postulation of nodes}.
Precisely, by using the notion of local positivity of line bundles on $X$, 
the machinery of Seshadri constants as in \cite{D}, \cite{EKL} and \cite{Ku} and 
the fundamental tool of the Kawamata-Viehweg vanishing theorem, 
we determine some upper-bounds for the number 
$\delta$, which are cubic polynomials in the integer $k$, such that 
if the $\delta$ nodes of a curve 
$C= V(s)$ are in very general position on $X$ (see Definition \ref{def:13}), then 
the point $[s]$ is regular for ${\V}_{\delta}({\E} \otimes \Oc_X(k))$
(see Theorem \ref{thm:14} and Corollary \ref{cor:14bis}). Furthermore, 
when the nodes of $C$ are assumed to be points in very general position
on an irreducible divisor of $X$ having at worst log-terminal 
singularities or to lie on a l.c.i and subcanonical curve in $X$, we determine upper-bounds
on $\delta$ which are, respectively, 
quadratic and linear polynomials in the integer $k$ implying the regularity of the point 
$[s]$ (see Theorems \ref{thm:19}, \ref{thm:21}, \ref{thm:34} and Corollaries
\ref{cor:20}, \ref{cor:35}).

We conclude the paper by focusing on the case $X= \Pt$ and by 
studying interesting geometric properties of space curves determined by elements 
in ${\V}_{\delta}({\E} \otimes \Oc_X(k))$.

The paper consists of six sections.
In Section \ref{S:1}, we recall some terminology and notation. Section 
\ref{S:2} contains fundamental definitions and technical details 
which are useful for our proofs. Section \ref{S:3} contains our main result
(see Theorem \ref{prop:3.fundamental}), which gives a cohomological 
description of the tangent space $T_{[s]}({\V}_{\delta} ({\F}))$. In Section
\ref{S:4}, we prove our uniform and effective result for the regularity of 
$ {\V}_{\delta} ({\F})$ at each point (see Theorem \ref{thm:9bis}). 
Section \ref{S:5} is devoted to the study of the regularity property of ${\V}_{\delta} ({\F})$
in terms of the postulation of nodes of the zero-loci of the
elements it parametrizes. We conclude with
Section \ref{S:6}, where we consider some geometric properties and
biliaison relation of space curves determined by elements in ${\V}_{\delta} ({\F})$.

{\it Acknowledgments:} Part of this paper was prepared during my permanence 
at the Department of Mathematics of the University of Illinois at Chicago 
(February - May 2001). Therefore, my deepest gratitude goes to L. Ein, not 
only for the organization of my visit, but mainly for all I have learnt from him, 
for having suggested me to approach this problem as well as for the 
valuable talks we had together. 
I am greatful to my collegues (and friends) G. Bini, N. Budur and T. de Fernex, 
for all they organized during my permanence and with whom I had some 
informal - but fundamental - talks. My warm thanks go 
to the Departments of Mathematics of the Universities of Illinois 
at Chicago and of Michigan at Ann Arbor for making my period of study and 
researches enjoyable, lively and productive. I am indebted to L. Chiantini, C. Ciliberto and 
the first referee of this paper for their remarks concerning Proposition \ref{prop:41}. 
My very special thanks go to GNSAGA-INdAM and to V. Barucci, A. F. Lopez and 
E. Sernesi for their confidence 
and their support during my period in U.S.A.

\section{Notation and Preliminaries}\label{S:1}
We work in the category of algebraic 
$\CC$-schemes. $Y$ is a \emph{$m$-fold} if it is a reduced, irreducible and non-singular scheme 
of finite type and of dimension $m$. 
If $m=1$, then $Y$ is a (smooth) {\em curve}; $m=2$ and $3$ are the cases of a 
(non-singular) {\em surface} and {\em threefold}, respectively. 
If $Z$ is a closed subscheme of a scheme $Y$, $\Ii_{Z/Y}$ 
denotes the \emph{ideal sheaf} of $Z$ in $Y$, ${\N}_{Z/Y}$ 
the {\em normal sheaf} of $Z$ in $Y$ whereas 
${\N}_{Z/Y}^{\vee} \cong {\Ii_{Z/Y}}/{\Ii_{Z/Y}^2}$ is the {\em 
conormal sheaf} of $Z$ in $Y$. As usual, $h^i(Y, \; -):=\text{dim} \; H^i(Y, \; -)$.

Given $Y$ a projective scheme, $\omega_Y$ denotes its dualizing sheaf. 
When $Y$ is a smooth variety, then $\omega_Y$ coincides with its canonical bundle and 
$K_Y$ denotes a canonical divisor s.t. $\omega_Y \cong \Oc_Y(K_Y)$. Furthermore, 
${\T}_Y$ denotes its tangent bundle whereas $\Omega^1_Y$ denotes its cotangent bundle. 

If $D$ is a reduced curve, $p_a(D)=h^1(\Oc_D)$ 
denotes its {\em arithmetic genus}, 
whereas $g(D)= p_g(D)$ denotes its \emph{geometric genus}, the
arithmetic genus of its normalization.

Consider $Y$ a projective $m$-fold. $Div(Y)$ denotes 
the set of (Cartier) divisors and $\sim$ the linear 
equivalence on $Div(Y)$, whereas $Pic(Y)$ denotes the Picard scheme of line bundles 
on $Y$. On the other hand, as in \cite{H1}, 
$F_1(Y)$ denotes the free abelian group generated by the set of all integral curves in $Y$. 
Denoted by $\cdot$ the \emph{intersection pairing} on $Y$ and 
by $\equiv$ the numerical equivalence on $Y$, we have
$$A^1(Y) = (Div(Y)/\equiv) \otimes_{\ZZ} \RR , \; {\rm and} \; 
A_1(Y) = (F_1(Y)/\equiv) \otimes_{\ZZ} \RR .$$Recall that 
an element $B \in Div(Y)$ is said to be {\em nef}, 
if $B \cdot D \geq 0$ for each irreducible curve $D$ on $Y$. 
A nef divisor $B$ is said to be {\em big} if $B^m>0$. 
By Kleiman's criterion (see, for example, \cite{H1}), a nef 
divisor $B$ is in the closure of the {\em ample divisor cone} 
$P^0(Y)$, which is the cone in $A^1(Y)$ generated by the ample divisors on $Y$. 
When $Y$ is a surface, in some literature, the ample divisor 
cone is also denoted by $N^+(Y)$ (see, for example, \cite{BPV} and \cite{Fr}).

Let $Y$ be a projective $m$-fold and $\E$ be a rank-$r$ vector bundle on 
$Y$; $c_i(\E)$ denotes the \emph{$i^{th}$-Chern class} of $\E$, 
$1 \leq i \leq r$. As in \cite{Ha} - Sect. II.7 - 
$\Pp_{Y} ({\E})$ denotes the {\it projective space bundle} on $Y$, 
defined as $Proj(Sym({\E}))$. There is a surjection 
$\pi^*(\E) \to \Oc_{\Pp_{Y} ({\E})}(1)$, where $\Oc_{\Pp_{Y} ({\E})}(1)$ is the 
{\it tautological line bundle} on $\Pp_{Y} ({\E})$ and where $\pi : \Pp_{Y} ({\E}) \to Y$ is 
the natural projection morphism. Recall that $\E$ is said to be an {\em ample} (resp. 
{\em nef}) vector bundle on $Y$ if $\Oc_{\Pp_{Y} ({\E})}(1)$ is an ample (resp. nef) 
line bundle on $\Pp_{Y} ({\E})$.

When $Y$ is a projective, normal variety of dimension $m$, 
the word {\em divisor} is used for Weil 
divisor, i.e. a formal linear combination of codimension-one subvarieties. A 
${\QQ}$-{\em divisor} ($\RR$-divisor, resp.) 
on $Y$ is a finite formal linear combination $D = \sum_i a_i D_i$ with 
rational (real, resp.) coefficients; when the coefficients are in $\ZZ$, $D$ is 
an {\em integral divisor}. $D$ is a $\QQ$-{\em Cartier} divisor if some multiple of $D$ is an
(integral) Cartier divisor (recall that, when $X$ is smooth, any $\QQ$-divisor is 
$\QQ$-Cartier). The {\em round-up} of $D$ and the 
{\em integral part} of $D$ are, respectively, the integral divisors
$\lceil D \rceil = \sum_i \lceil a_i \rceil D_i$ and 
$[D] = \sum_i [a_i] D_i$ where, as usual, for $x \in \QQ$ one denotes by $\lceil x \rceil$ 
the least integer greater than or equal to $x$ and by $[x]$ the greatest integer 
smaller than or equal to $x$. The {\em fractional part} of $D$ is 
$\{D\} = D - [D]$. Since there is a $\QQ$-valued intersection theory for 
$\QQ$-Cartier $\QQ$-divisors, one can extend the notion of ampleness and nefness 
to $\QQ$-divisors. Similarly, $D$ is big if $nD$ is integral and big, for some positive $n$. 
$D=\sum_i a_i D_i$ has {\em simple normal crossings} if each $D_i$ is smooth 
and if $D$ is defined in a neighborhood of any point by an equation in local 
analytic coordinates of the type $z_1 \cdots z_k = 0$, with $k \leq m$. A 
{\em boundary divisor} $\Delta$ is an effective divisor whose support 
has simple normal crossings and such that $[\Delta] = 0$.

If $Y$ is a projective normal variety and if $D$ is a $\QQ$-divisor on $Y$, a 
{\em log-resolution} of the pair $(Y, D)$ is a proper birational mapping 
$\mu : Y' \to Y,$ where $Y'$ is smooth and such that 
the divisor $\mu^*(D) + Exc(\mu)$ has simple normal crossing support 
(here $Exc(\mu)$ denotes the sum of the $\mu$-exceptional divisors). If $(Y,D)$ is such that 
$K_Y + D$ is $\QQ$-Cartier and if $\mu$ is a log-resolution of the pair, then, 
$$K_{Y'/Y} - \mu^*(D) := K_{Y'} - \mu^*(K_Y + D) \equiv \sum_i a_i E_i,$$where the $E_i$'s are 
distinct irreducible divisors (not necessarily all $\mu$-exceptionals) and 
where the coefficients $a_i$'s are called the {\em discrepancies}. The pair $(Y,D)$ is 
called {\em log-terminal} (resp. {\em Kawamata log-terminal}) if 
$a_i > -1$ for each $E_i$ $\mu$-exceptional (resp. for each $E_i$). Since 
$K_{Y'/Y}$ is always $\mu$-exceptional, 
$Y$ is said to have {\em at worst log-terminal singularities} if the 
pair $(Y,0)$ is log-terminal.

To conclude, we also recall one of the most 
important vanishing theorem for $\QQ$-divisors - the Kawamata-Viehweg theorem - 
which will be frequently used in the sequel. 

\noindent
{\bf Theorem} (Kawamata-Viehweg, see, for example, \cite{Miy}, page 146) 
{\em Let $Y$ be a smooth, projective variety of dimensione $m$, let $D$ 
be a big and nef $\QQ$-divisor, whose fractional part has simple normal 
crossing support. Then
$$H^i(Y, {\Oc}_Y(K_Y + \lceil D \rceil ) = (0), \; {\rm for} \; i >0.$$When 
$Y$ is a surface, the same conclusion holds even without the hypothesis on 
the fractional part of $D$.}

\section{Basic definitions and fundamental properties}\label{S:2}

In this section we introduce some fundamental definitions 
and remarks concerning the study of families of curves on smooth projective 
threefolds. For generalities, the reader is referred to \cite{Sz}, Chapter IV.

\begin{definition}\label{def:koszul}
Let $X $ be a smooth projective threefold and let $\F$ be a rank-two vector 
bundle on $X$. Let $s$ be a global section of $\F$. The {\em zero-locus} of $s$, 
denoted by $V(s)$, is the closed subscheme of $X$ defined by the exact sequence
$${\F}^{\vee} \stackrel{s^{\vee}}{\longrightarrow} {\Oc}_X \to {\Oc}_{V(s)} \to 0,$$where 
$s^{\vee}$ is the dual map of the section $s$.  
\end{definition}

\noindent
If $codim_X(V(s)) = 2$, then $Ker(s^{\vee}) = L^{\vee}$ is a line bundle on $X$ such 
that $c_1({\F}) = \bigwedge^2({\F}) \cong L$. This yields the {\em Koszul sequence} 
of $({\F}, s)$:
\begin{equation}\label{eq:koszul}
0 \to {\Oc}_X \to {\F} \to {\Ii}_{V(s)} \otimes L \to 0.
\end{equation}

\begin{remark}\label{rem:pic}
\normalfont{
When $Pic(X) \cong \ZZ$ (e.g $X= \Pt$ or $X$ either a prime Fano or a complete intersection 
Calabi Yau threefold) one can use this isomorphism to identify line bundles on $X$ with integers. 
In particular, if $A$ denotes the ample generator class of $Pic(X)$ over $\ZZ$ and 
if $\F$ is a rank-two vector bundle on $X$ such that $c_1({\F}) = n A$, we can also 
write $c_1({\F}) = n$ with no ambiguity.
}
\end{remark}

\noindent
We recall well-known results concerning the correspondence 
between curves and global sections of vector bundles on a smooth projective threefold.

\begin{theorem}\label{thm:serre}(Serre)
Let $X$ be a smooth projective threefold. A curve $D \subset X$ occurs as the zero-locus 
of a global section of a rank-two vector bundle $\F$ on $X$ if and only if $D$ is locally complete 
intersection and its dualizing sheaf $\omega_D$ is isomorphic to the restriction to $D$ of 
$\omega_X \otimes M$, for some line bundle $M$ on $X$ such that 
\begin{equation}\label{eq:condition}
h^1(X, M^{\vee}) = h^2(X, M^{\vee})=0.
\end{equation}Furthermore, such a curve $D$ is a complete intersection in $X$ iff $\F$ 
splits.
\end{theorem}

\noindent
\begin{proposition}\label{prop:utile}
Let $X$ be a smooth projective threefold and let $\F$ be a rank-two vector 
bundle on $X$. If $\F$ is globally generated, then the zero-locus of a general section of $\F$, if 
not empty, is non-singular of codimension two.
\end{proposition} 
\begin{proof}This is a particular case of standard results - extending Bertini's theorem - 
concerning degeneracy-loci of generic morphisms between vector bundles on $X$.
\end{proof}

From what recalled above, if $X$ is a smooth projective threefold and if $\F$ is a 
globally generated rank-two vector bundle on $X$, it is not restrictive if from now on 
we assume 
that the zero-locus of the 
general section of $\F$ is a smooth, irreducible curve $D$ in 
$X$. By the Koszul sequence (\ref{eq:koszul}), 
we find the geometric genus of $D$ in terms of the Chern classes of $\F$ and of the 
invariants of $X$. Precisely 
\begin{equation}\label{eq:numeriX}
2g(D) -2 = 2p_a(D) -2= deg(c_1({\F}) \otimes \omega_X \otimes 
{\Oc}_D). 
\end{equation}This integer is easily 
computable when, for example, $X$ is a general complete intersection threefold. 
In particular, when $X= \Pt$, by Remark \ref{rem:pic}, if 
we put $c_i = c_i({\F}) \in \ZZ$, we have
\begin{equation}\label{eq:numeriP3}
deg(D) = c_2 \; {\rm and} \; g(D)= p_g(D) = \frac{1}{2} (c_2 (c_1 -4))+ 1, 
\end{equation}i.e. $D$ is subcanonical of level $(c_1 -4)$. 

Take now ${\Pp}(H^0({\F}))$; from our assumptions on $\F$, 
the general point of this projective space 
parametrizes a global section whose zero-locus is a 
smooth, irreducible curve in $X$. 
Given a positive integer $\delta$, one can consider the subset
\begin{equation}\label{eq:severi var}
\begin{aligned}
{\V}_{\delta}({\F}) := & \{[s] \in {\Pp}(H^0({\F})) \; | \; C_s := V(s) \subset X \; {\rm is \; irreducible} \\
 & {\rm  with \; only} \;  \delta \; {\rm nodes \; as \; singularities} \}; 
\end{aligned}
\end{equation}therefore, any element of ${\V}_{\delta}({\F})$ determines a curve in $X$ 
whose arithmetic genus $p_a(C_s)$ is given by 
(\ref{eq:numeriX}) and whose geometric genus is $g= p_a(C_s) - \delta$.

${\V}_{\delta}({\F})$ is a locally closed 
subscheme of the projective space ${\Pp}(H^0({\F}))$ and it 
is usually called the {\em Severi variety} of global sections of $\F$ 
whose zero-loci are irreducible, 
$\delta$-nodal curves in $X$ (see \cite{BC}); this is because such schemes are 
the natural generalization of the (classical) Severi varieties of irreducible and $\delta$-nodal curves 
in linear systems on smooth, projective surfaces (see \cite{CC}, \cite{CH},
\cite{CS}, \cite{F1}, \cite{GLS}, \cite{Harris}, \cite{Ran}, \cite{S} and \cite{Sev}, just to 
mention a few).  

\begin{proposition}\label{prop:0}
Let $X $ be a smooth projective threefold, $\F$ be a 
globally generated rank-two vector bundle on $X$ and $\delta$ be a positive integer. 
Then
\[expdim ({\V}_{\delta}({\F})) = \begin{cases}
				 h^0(X, {\F}) - 1 - \delta, \; \; {\rm if} \; \delta \leq 
                                 h^0(X, {\F}) - 1= dim({\Pp}(H^0({\F}))),\\ 
 				 -1, \;\;\;\;\;\;\;\;\;\;\;\;\;\;\;\;\;\;\;\;\;\;\;\;{\rm if} 
                                 \; \delta \geq 
                                 h^0(X, {\F}).\\ 
				\end{cases} \]

\end{proposition}
\begin{proof}
If ${\V}_{\delta}({\F}) = \emptyset$, then $dim ({\V}_{\delta}({\F})) = -1$. On the other 
hand, when ${\V}_{\delta}({\F}) \neq \emptyset$, consider the set 
$${\mathcal U}_{\delta} \subset {\Pp}(H^0({\F})) \times (X^{\delta} \setminus 
\bigcup_{1 \leq i \neq j \leq \delta}\Delta_{i,j}),$$where 
$X^{\delta}$ is the $\delta$-Cartesian product of $X$, $\Delta_{i,j}$ are the 
diagonals in $X^{\delta}$ and where 
$${\mathcal U}_{\delta} := \{([s]; \; p_1, \ldots , p_{\delta}) \; | \;
C = V(s) \subset X \; {\rm is \; an} \; {\rm irreducible \; curve \; 
with \, only \; nodes \, at \; the} \; p_i's\}.$$Since $X$ is smooth, for an arbitrary 
$p \in X$ we consider $U = U_p$ an affine open subscheme of $X$ containing $p$, with 
$(x_1^{(p)},x_2^{(p)},x_3^{(p)})$ local coordinates in $U_p$, 
such that $s |_{U_p} = (f_1^{(p)}, f_2^{(p)})$, where $f_i^{(p)} \in {\Oc}_X(U_p)$. Define 
the closed subscheme
\begin{displaymath}
\begin{array}{cl}
{\mathcal K}_{\delta}:= &\{([s]; \; p_1, \ldots , p_{\delta}) \in 
{\Pp}(H^0({\F})) \times (X^{\delta} \setminus 
\bigcup_{1 \leq i \neq j \leq \delta}\Delta_{i,j}) \; | \; s(p_i) = 0 \;  {\rm and}\\
 & \;rank(J(s)(p_i)) \leq 1, \; 1 \leq i \leq \delta \},
\end{array}
\end{displaymath}where, for an arbitrary $p \in X$, 
$J(s)(p)$ is the Jacobian matrix of $s$ at the point $p$. By definition, 
\begin{displaymath}
\begin{array}{cl}
{\mathcal K}_{\delta}:= & \{([s]; p_1, \ldots , p_{\delta}) \; | \; s(p_i) = 
((\frac{\partial}{\partial x_1^{(p_i)}} \wedge \frac{\partial}{\partial x_2^{(p_i)}})(s)) 
(p_i) = \\
& \; ((\frac{\partial}{\partial x_1^{(p_i)}} \wedge \frac{\partial}{\partial x_3^{(p_i)}})(s)) (p_i) 
= ((\frac{\partial}{\partial x_2^{(p_i)}} \wedge \frac{\partial}{\partial x_3^{(p_i)}})(s)) (p_i) = 0, \; 
1 \leq i \leq \delta \}.
\end{array}
\end{displaymath}Since ${\mathcal U}_{\delta}$ is contained in 
${\mathcal K}_{\delta}$ as an open dense subscheme and since 
${\mathcal K}_{\delta}$ is cut out by at most $4 \delta $ independent 
equations, then
\begin{displaymath}
\begin{array}{ll}
dim({\mathcal U}_{\delta}) =  dim({\mathcal K}_{\delta}) & \geq  dim({\Pp}(H^0({\F})) 
\times (X^{\delta} \setminus 
\bigcup_{1 \leq i \neq j \leq \delta}\Delta_{i,j})) - 4 \delta = \\
& = h^0(X, {\F}) - 1 + 3 \delta - 4 \delta =  h^0(X, {\F}) - 1 - \delta .
\end{array}
\end{displaymath}Denoting by $\pi_1$ the restriction to 
${\mathcal U}_{\delta}$ of the projection onto the 
first factor of the product $ {\Pp}(H^0({\F})) \times  X^{\delta}$, 
we have $\pi_1({\mathcal U}_{\delta}) = {\V}_{\delta}({\F})$. We conclude by 
observing that $\pi_1$ is finite onto its image. 
\end{proof}   

\noindent
{\bf Assumption}: From now on, we shall use Notation \ref{not:1}. Moreover,
given $X$ and $\F$ as in Proposition \ref{prop:0},
we shall always assume ${\V}_{\delta}({\F}) \neq \emptyset$ and 
$\delta \leq {\rm min} \{ h^0(X, {\F}) - 1, \; p_a(C) \}$ (the latter 
is because we want $C=V(s)$ to be irreducible).

By Proposition \ref{prop:0}, we can state the following fundamental definition.

\begin{definition}\label{def:0}
Let $[s] \in {\V}_{\delta}({\F})$, with 
$\delta \leq {\rm min} \{ h^0(X, {\F}) - 1, \; p_a(C) \}$.
Then $[s]$ is said to be a {\em regular point} of ${\V}_{\delta}({\F})$ if:
\begin{itemize}
\item[(i)] $[s] \in {\V}_{\delta}({\F})$ is a smooth point, and
\item[(ii)] $dim_{[s]}({\V}_{\delta}({\F})) = expdim ({\V}_{\delta}({\F}))=
h^0(X, {\F}) - 1 - \delta$, i.e. $$dim_{[s]}({\V}_{\delta}({\F})) =
dim({\Pp}(H^0({\F}))) - \delta.$$
\end{itemize}
\end{definition}The goal of the next section is to
present a cohomological description
of the tangent space $T_{[s]}({\V}_{\delta}({\F}))$ which will translate
the regularity property of the point $[s] \in {\V}_{\delta}({\F})$ into the
surjectivity of some maps among vector spaces of sections of suitable sheaves
on the threefold $X$.

\section{Description of the tangent space $T_{[s]}({\V}_{\delta}({\F}))$ and
regularity}\label{S:3}

As before, let $X$ be a smooth projective threefold and let $\F$ be a globally generated
rank-two vector bundle on $X$. Let $\delta$ be a positive integer and consider
$[s] \in {\V}_{\delta}({\F})$. From now on in this section,
let $C$ be the curve in $X$ which is the zero-locus of the given $s$ and let
$\Sigma$ denote its set of $\delta$ nodes.

Since $C$
is local complete intersection in $X$, its normal sheaf
${\N}_{C/X}$ is a rank-two vector-bundle (see \cite{Ha}). Precisely,
${\N}_{C/X} \cong {\F}|_C$. Let $T^1_C$ be the {\em first cotangent sheaf} of $C$, i.e.
$T^1_C \cong {\mathcal Ext}^1(\Omega^1_C, {\Oc}_C)$,
where $\Omega_C^1$ is the sheaf of K$\ddot{a}$hler differentials of the nodal
curve $C$ (for details, see \cite{LS}). We have the exact sequence
\begin{equation}\label{eq:T1}
0 \to {\N}'_C \to {\N}_{C/X} \stackrel{\gamma}{\to} T^1_C \to 0, 
\end{equation}where ${\N}'_C$ is defined as the kernel of the natural surjection 
$\gamma$ (see, for example, \cite{S}). Since nodal points are planar 
singularities, one has 
\begin{itemize}
\item[i)] $T^1_{C,p} = 0$ and ${\N}'_{C,p} \cong  {\N}_{C/X,p} \cong {\Oc}_{C,p}^{\oplus \;2}$, when 
$p \in C$ is a smooth point,
\item[ii)] $T^1_{C,p}\cong \CC$ and ${\N}'_{C,p} \cong  (\underline{m}_p {\Oc}_{C,p}) 
\oplus {\Oc}_{C,p}$, when $p$ is a node of $C$ ($\underline{m}_p$ denotes the maximal ideal 
at the point $p$). 
\end{itemize}Therefore, $T^1_C$ is a sky-scraper sheaf supported on $\Sigma$, such that
$T^1_C \cong \bigoplus_{i=1}^{\delta} {\CC}_{(i)}$. 

By using (\ref{eq:T1}), the goal of this section is to construct a subsheaf 
${\F}^{\Sigma} \subset {\F}$ fitting in 
the following exact diagram:
\begin{equation}\label{eq:(*1)} 
\begin{aligned}
\begin{array}{rcccclr}
 &0 & & 0 & & & \\
 & \downarrow & & \downarrow & & & \\
0 \to   &{\Ii}_{C/X} \otimes {\F}  & \stackrel{\cong}{\to} & {\Ii}_{C/X}\otimes {\F} &\to & 0
&  \\ 
 & \downarrow & & \downarrow & & \downarrow & \\
 0 \to  & {\F}^{\Sigma} & \to & {\F} & \to  & T^1_C & \to 0 \\ 
 & \downarrow & & \downarrow & & \downarrow^{\cong}& \\
0 \to & {\N}'_C & \to &  {\F} |_{C} & \to & T^1_C & \to 0  \\ 
 &\downarrow & & \downarrow & &\downarrow & \\
 & 0 & & 0 & & 0 & .
\end{array}
\end{aligned}
\end{equation} Observe that, from the commutativity of diagram \eqref{eq:(*1)}, 
$H^0(X, {\F}^{\Sigma})/ < s >$ 
parametrizes the first-order deformations of the section $s$ in $H^0(X, {\F}) $ which are 
equisingular; 
indeed, these are exactly the global sections of ${\F}$ which go to zero at $\Sigma$ in the 
composition 
\begin{equation}\label{eq:compo}
{\F} \to {\F}|_C \to T^1_C \cong {\Oc}_{\Sigma} \to 0.
\end{equation}

In the following result, which is the core of the entire paper,
we construct the sheaf ${\F}^{\Sigma}$ by using some projective space-bundle arguments.

\begin{theorem}\label{prop:3.fundamental}
Let $X $ be a smooth projective threefold. Let $\F$ be a globally generated rank-two vector
bundle on $X$ and let $\delta$ be a positive integer. As in \eqref{eq:severi var}, let
\[\begin{aligned}
{\V}_{\delta}({\F}) = & \{[s] \in {\Pp}(H^0({\F})) \; | \; C_s := V(s)
\subset X \; {\rm is \; irreducible \;} \\
 & {\rm  with \; only} \;  \delta \; {\rm nodes \; as \; singularities} \}.
\end{aligned}\]Fix $[s] \in {\V}_{\delta}({\F})$ and
let $C=V(s) \subset X$. Denote by $\Sigma$ the set of nodes of
$C$. Let$${\mathcal P} := {\Pp}_{X}({\F}) \stackrel{\pi}{\longrightarrow} X$$be
the
projective space bundle together with its natural projection $\pi$ on $X$ and denote
by ${\Oc}_{\mathcal P}(1)$ its tautological line bundle. Let
$$\Sigma^1 : = {\Pp}_X(T^1_{C}) \subset {\Pmc}$$denote the
zero-dimensional subscheme of $\Pmc$ of length $\delta$, determined by the
surjection \eqref{eq:compo}. Denote by $\Ii_{\Sigma^1/{\Pmc}}$ the ideal sheaf of 
$\Sigma^1$ in $\Pmc$. Then
\begin{itemize}
\item[(i)] $\Sigma^1$ is a set of $\delta$ rational double points for the
divisor $D_s \in |{\Oc}_{\mathcal P}(1)|$, corresponding to the given section
$s \in H^0(X, \F)$, and
\item[(ii)] the subsheaf of $\F$, defined by
\begin{equation}\label{eq:svolta}
{\F}^{\Sigma} := \pi_* ({\Ii}_{\Sigma^1/{\Pmc}} \otimes {\Oc}_{\mathcal P}(1)),
\end{equation}is such
that its global sections (modulo the
one dimensional subspace $<s>$) parametrize first-order deformations of
$s \in H^0(X, \F)$ which are equisingular.

\noindent
In particular, we have
\begin{equation}\label{eq:reg}
\frac{H^0(X, {\F}^{\Sigma})}{< s >} \cong T_{[s]} ({\V}_{\delta}({\F})) \subset T_{[s]}
({\Pp}(H^0({\F}))) \cong \frac{H^0(X, {\F})}{< s >}.
\end{equation}
\end{itemize}
\end{theorem}
\begin{proof}
To naturally define the sheaf ${\F}^{\Sigma}$ and the diagram \eqref{eq:(*1)}, we
consider the smooth, projective fourfold
$${\mathcal P} := {\Pp}_{X}({\F}) \stackrel{\pi}{\longrightarrow} X,$$together
with its tautological line bundle ${\Oc}_{\mathcal P}(1)$ such that
$\pi_*({\Oc}_{\mathcal P}(1)) \cong {\F}$.
From$$0 \to {\Oc}_X \stackrel{\cdot s}{\to} {\F},$$we also have
$$0 \to {\Oc}_{\Pmc} \stackrel{\cdot s}{\to} {\Oc}_{\Pmc}(1).$$Therefore, the nodal curve
$C \subset X$ corresponds to a divisor $D_s \in |{\Oc}_{\Pmc}(1) |$ on the
fourfold $\Pmc$. Take also
$${\FF} := {\Pp}^1_C = {\mathcal Proj} ({\Oc}_C[\xi_0, \xi_1]) \stackrel{\pi_1}{\to} C$$which
is a ruled surface in $\Pmc$. We want to study some geometric properties of $D_s$
and of $\FF$. Let $p \in \Sigma = Sing(C)$. Take $U_p \subset X$
an affine open set containing $p$, where the vector bundle ${\F}$ trivializes.
Choose local coordinates $\underline{x} = (x_1, x_2, x_3)$
on $U_p \cong {\CC}^3$ such
that $\underline{x}(p) = (0,0,0)$ and such that the global section $s$ is
$$s|_{U_p} = (x_1x_2, \; x_3).$$Then
$${\Oc}_C(U_p) \cong  {\CC}[x_1,x_2,x_3]/(x_1x_2, \; x_3),$$and
$${\Oc}_{{\FF}}(\pi_1^{-1}(U_p)) \cong  {\CC}[x_1,x_2,x_3, \frac{\xi_1}{\xi_0}]/(x_1x_2,
\; x_3).$$Therefore, the surface
$\FF$ is singular along the lines in ${\La} = \bigcup_{i=1}^{\delta} L_i = \pi_1^{-1} (\Sigma)
= \pi^{-1} (\Sigma)$. For what concerns $D_s \in |{\Oc}_{\mathcal P}(1)|$, since
$U_p$ trivializes $\F$, then ${\Pmc} \; |_{U_p} \cong U_p \times {\Pp}^1.$ Taking homogeneous
coordinates $[u,v] \in {\Pp}^1$,we have ${\Oc}_{\Pmc}(U_p) \cong {\CC}[x_1, x_2, x_3, u, v].$
Thus,
$$ {\Oc}_{D_s}(U_p) \cong {\CC}[x_1, x_2, x_3, u, v]/(u x_1 x_2 + v x_3).$$By standard
computations,
we see that $D_s$ has a rational double point along the line $\pi^{-1} (p) = \pi_1^{-1}(p)$ which
belongs to the singular locus of ${\FF} \subset \Pmc$.

Globally speaking, by using \eqref{eq:compo}, we can state that
the divisor $D_s \subset \Pmc$ is singular along
the locus $$\Sigma^1 : = {\Pp}_X(T^1_C) \subset {\Pmc} = {\Pp}_X ({\F}),$$where $\Sigma^1 \cong
\Sigma$ is a set of $\delta$ rational double points for $D_s$, each line of
${\La} = \pi^{-1}(\Sigma)$ containing only one of such $\delta$ points. Since
$ \Sigma^1 \subset \Pmc$ is a closed immersion, we have the natural exact sequence
\begin{equation}\label{eq:idealP}
0 \to {\Ii}_{\Sigma^1/{\Pmc}} \otimes {\Oc}_{\mathcal P}(1) \to {\Oc}_{\mathcal P}(1)
\to {\Oc}_{\Sigma^1} \to 0,
\end{equation}which is defined by restricting ${\Oc}_{\mathcal P}(1) $ to $\Sigma^1$. By the
definition of tautological line bundle, we have:
\begin{displaymath}
\begin{array}{rcccccc}
 & & & \pi^*({\F}) & \to  & \pi^*({\Oc}_{\Sigma})& \to 0 \\
 &  & & \downarrow & & \downarrow& \\
0 \to & {\Ii}_{\Sigma^1/{\Pmc}} \otimes {\Oc}_{\mathcal P}(1) & \to &  {\Oc}_{\mathcal P}(1)  &
\to & {\Oc}_{\Sigma^1} & \to 0 \\
 &\downarrow & & \downarrow & &\downarrow & \\
 & 0 & & 0 & & 0 &.
\end{array}
\end{displaymath}Since $$\pi_*({\Oc}_{\mathcal P}(1)) \cong {\F}, \;
\pi_*({\Oc}_{\Sigma^1}) = \pi_*({\Oc}_{\pi^{-1}(\Sigma^1)}) = \pi_*(\pi^*({\Oc}_{\Sigma}))
\cong {\Oc}_{\Sigma}$$and since we have
${\F} \to \!\!\! \to {\Oc}_{\Sigma}$, by applying $\pi_*$ to the exact sequence
(\ref{eq:idealP}), we get ${\mathcal R}^1 \pi_*({\Ii}_{\Sigma^1/{\Pmc}}
\otimes {\Oc}_{\mathcal P}(1)) = 0$. Thus, we define
$${\F}^{\Sigma} := \pi_* ({\Ii}_{\Sigma^1/{\Pmc}} \otimes {\Oc}_{\mathcal P}(1)),$$which
gives \eqref{eq:svolta}, so that
\begin{equation}\label{eq:idealX}
0 \to   {\F}^{\Sigma}  \to {\F} \to   {\Oc}_{\Sigma} \to 0,
\end{equation}as well as diagram \eqref{eq:(*1)}, holds.
\end{proof}

We remark that \eqref{eq:reg} gives a completely general characterization of the tangent 
space $T_{[s]} ({\V}_{\delta}({\F}))$ on $X$. Furthermore, we have:

\begin{coroll}\label{rem:reg}
With assumptions as in Theorem \ref{prop:3.fundamental}, from \eqref{eq:idealX} we get
\begin{equation}\label{eq:regbis}
[s] \in {\V}_{\delta}({\F}) \; {\rm is \; regular} \; \Leftrightarrow
H^0(X , {\F}) \stackrel{\alpha_X}{\to \!\!\! \to} H^0(X, {\Oc}_{\Sigma}) \Leftrightarrow
H^0({\Pmc} , {\Oc}_{\mathcal P}(1)) \stackrel{\alpha_{\Pmc}}{\to \!\!\! \to} H^0({\Pmc},
{\Oc}_{\Sigma^1}).
\end{equation}
\end{coroll}
\begin{proof}
It follows from Proposition \ref{prop:0} and from Theorem \ref{prop:3.fundamental}.
\end{proof}

Note that, on the one hand, the map $\alpha_X$ in \eqref{eq:regbis} is not defined by restricting 
the global sections of $\F$ to $\Sigma$ because (\ref{eq:idealX}) - i.e.
the second row of diagram \eqref{eq:(*1)} -
does not coincide with the restriction sequence
$$0 \to   {\Ii}_{\Sigma/X} \otimes {\F}  \to {\F} \to   {\F} |_{\Sigma} \to 0;$$precisely, 
we have
\begin{equation}\label{eq:circledast}
\begin{array}{rcccccr}
 &0 & & 0 & & 0 & \\
 & \downarrow & & \downarrow & & \downarrow & \\ 
0 \to   & {\Ii}_{\Sigma/X} \otimes {\F}  & \to & {\F}^{\Sigma} & \to & {\Oc}_{\Sigma}
 & \to 0 \\
 & \downarrow^{\cong} & & \downarrow & & \downarrow & \\ 
 0 \to  & {\Ii}_{\Sigma/X} \otimes {\F} & \to & {\F} & \to  
& {\F}|_{\Sigma} \cong 
{\Oc}_{\Sigma}^{\oplus 2} & \to 0 \\ 
 & \downarrow & & \downarrow & & \downarrow & \\ 
& 0 & \to &  {\Oc}_{\Sigma} & \to &{\Oc}_{\Sigma} & \to 0 ;\\ 
 & & & \downarrow & &\downarrow & \\ 
 &  & & 0 & & 0 & 
\end{array}
\end{equation}

On 
the other hand, the exact sequence (\ref{eq:idealP}) on the fourfold $\Pmc$ is equivalent 
to (\ref{eq:idealX}), by the Leray isomorphisms, but it is more naturally defined
by restricting the line bundle ${\Oc}_{\mathcal P}(1)$ to $\Sigma^1$. Therefore, the map 
$\alpha_{\Pmc}$ in (\ref{eq:regbis}) is a classical restriction map.

To better understand the map $\alpha_X$, we also want to give a local description of 
(\ref{eq:idealX}).

\vskip 10pt

\noindent
\underline{Local description}: let $p \in Sing(C) = \Sigma$ and take, as before, 
$U_p \subset X$ an affine open set containing $p$, where the vector bundle ${\F}$ is trivial. 
Take local 
coordinates $\underline{x} = (x_1, x_2, x_3)$ 
on $U_p \cong {\CC}^3$ such 
that $\underline{x}(p) = (0,0,0)$ and such that the global section $s$, whose zero-locus 
is $C$, is $s|_{U_p} = (x_1x_2, \; x_3).$ Since $C = V(x_1x_2, \; x_3) \subset 
Spec({\CC}[x_1, x_2, x_3])  \cong U_p$, around the node 
$ \underline{x}(p)= \underline{0} $ the map
$$(**) \;\;\;\;  {\T}_{{\CC}^3}|_C \stackrel{J(s)}{\longrightarrow} {\N}_{C/{\CC}^3} \to T^1_C $$is given 
by  
\[J(s) := \left( \begin{array}{ccc}
			\frac{\partial f_1}{\partial x_1} & \frac{\partial f_1}{\partial x_2} & 
                         \frac{\partial f_1}{\partial x_3}\\
                         \frac{\partial f_2}{\partial x_1} & \frac{\partial f_2}{\partial x_2} & 
                         \frac{\partial f_2}{\partial x_3}
		 	\end{array}
			\right) = 
                    \left( \begin{array}{ccc}
			x_2 & x_1 & 
                         0 \\
                       0  & 0  & 1
		 	\end{array}
			\right). \]From the fact that 
${\rm rank} (J(s)|_{\underline{0}}) = 1$, it follows that 
$coker(J(s)|_{\underline{0}}) \cong \CC$. Let $s(x_1, x_2, x_3 ) = s|_{U_p} = 
 (x_1x_2, \; x_3)$; if $\sigma (x_1, x_2, x_3 )$ is 
a section of ${\F}^{\Sigma}$ over $U_p$ then, by definition, 
$$s_{\epsilon}(x_1, x_2, x_3 ): = s(x_1, x_2, x_3 ) + \epsilon \; \sigma (x_1, x_2, x_3 ) 
$$is a first-order deformation of $s$ which determines equisingular zero-loci. Then, by $(**)$,
$$\sigma (x_1, x_2, x_3 ) = J(s) \underline{u},$$where
$\underline{u} = \underline{u}(x_1, x_2, x_3) = (u_1(x_1, x_2, x_3),
u_2(x_1, x_2, x_3), u_3 (x_1, x_2, x_3))$. To see this, consider
$$(***) \;\;\;\;\; s(\underline{x} + \epsilon \underline{u}) \; |_{\underline{0}} \equiv
(s(\underline{x}) + \epsilon J(s) \underline{u}) \; |_{\underline{0}} \pmod{\epsilon^2}
= \left( \begin{array}{c}
		x_1x_2 \\
		x_3
	\end{array}
            \right)
+ \epsilon  \left( \begin{array}{c}
		x_2 u_1 + x_1 u_2 \\
		u_3
	\end{array}
        \right) \pmod{\epsilon^2};$$thus
$$s_{\epsilon}(\underline{x} )  \equiv s(\underline{x} + \epsilon \underline{u})
\pmod{\epsilon^2}.$$Moreover, since $U_p$ is a trivializing open subset for $\F$,
we have that $(**)$ becomes
\begin{displaymath}
\begin{array}{cccr}
 {\Oc}_C^{\oplus 3} & \stackrel{J(s)}{\longrightarrow}& {\Oc}_C^{\oplus 2} & \to T^1_C \\
(e_1, e_2, e_3) &  \to & (e'_1, e'_2) & .
\end{array}
\end{displaymath}Since $Im(J(s)) = < x_2 e'_1, x_1 e'_1, e'_2>$, then
$e'_2$ goes to zero in $T^1_C$ so the deformations in $(***)$ are actually equisingular.

Observe that, by (\ref{eq:regbis}) and by using some results in \cite{deC},
one can immediately determine some conditions for the regularity of
$[s]\in {\V}_{\delta}({\F})$. Indeed, if $\F$ is a globally generated
rank-two vector bundle on $X$ which generates the $0$-jets at 
$\Sigma = \{p_1, \; \ldots, \; p_{\delta}\}$ (equiv. which separates 
the points of $\Sigma$), then by definition we have
$$H^0(X, {\F}) \to \!\!\! \to \bigoplus_{i=1}^{\delta}{\CC}^2_{(i)}.$$This implies 
the regularity conditions (\ref{eq:regbis}), as it immediately follows 
by considering the last two columns of \eqref{eq:circledast}. 
In such a case $[s] \in {\V}_{\delta}({\F})$,
such that $C = V(s)$ and $\Sigma = Sing(C)$, is therefore a regular point.
From Proposition 3.1 in \cite{deC}, we deduce:

\begin{proposition}\label{prop:1.bis}
Let $X$ be a smooth, projective threefold. Let $\G$ be a nef rank-two vector bundle and
$L$ be a big and nef line bundle on $X$. Consider the rank-two vector bundle
$${\F} := {\G} \otimes \omega_X \otimes det({\G}) \otimes {\Oc}_X(L)$$whose general
section is assumed to be a smooth curve in $X$. Take $\delta$ be a positive integer and consider
$[s] \in {\V}_{\delta}({\F})$. Let $\Sigma = \{ p_1, \; \ldots , p_{\delta} \}$ be the set
of nodes of $C= V(s)$. Let ${\epsilon}(L, p)$ denote the Seshadri constant of the line bundle
$L$ at the point $p$ (for precise definition see Remark \ref{rem:3} or Definition
\ref{def:12}).
Assume either
$${\epsilon}(L, p_j) > 3 {\delta} , \; \forall \; p_j \in \Sigma$$or
$$L^3 > ({\epsilon}(L, p_j))^3 \; {\rm and} \; {\epsilon}(L, p_j) \geq 3 {\delta},\;
\forall \; p_j \in \Sigma .$$Then, the global
sections of $\F$ separate $\Sigma$. In particular, ${\V}_{\delta}({\F})$ is regular at $[s]$.
\end{proposition}
\begin{proof}
One applies the effective non-vanishing Theorem 2.2 in \cite{deC} taking
$L_j = \frac{1}{\delta}L$.
\end{proof}

\begin{remark}\label{rem:3}
{\normalfont
Take e.g. $X \subset \PR$ a smooth threefold, whose hyperplane
section is denoted by $H$. Consider the line bundle
${\Oc}_X(k H)$, where $k$ is a positive integer, and
take $\Sigma = \{p_1, \ldots, p_{\delta} \} \subset X$. Denote by
$\mu_j$ the blowing-up of $X$ at the point $p_j$. Then, by definition,
$$\epsilon({\Oc}_X(kH) , p_j):=
{\rm Sup} \{ \epsilon \in \RR_{\geq 0} | \; \mu_j^*(kH) - \epsilon E_j \; {\rm
is \; a \; nef \; {\RR}-divisor \; on} \; Bl_{p_j}(X) \},$$
where $E_j$ denotes the $\mu_j$-exceptional divisor. Equivalently,
$$\epsilon({\Oc}_X(kH) , p_j):= {\it Inf}_{\Gamma \subset X}
\{ \frac{kH \cdot \Gamma}{mult_{p_j}(\Gamma)} \}$$where the infimum is
taken over all reduced and irreducible curves
$\Gamma \subset X$ passing through $p_j$. Since $H$ is very ample on $X$, then the
curve $\Gamma$ - as a curve in $\PR$ - is such that
$deg(\Gamma) \geq mult_q(\Gamma)$, for each $q \in \Gamma$. Therefore, the numerical
conditions in Proposition \ref{prop:1.bis} give $$\delta \leq \frac{k}{3},$$
which is a linear bound on the admissible number of nodes of $C=V(s)$ in order to have that
$[s] \in {\V}_{\delta} ({\G} \otimes \omega_X \otimes det({\G}) \otimes {\Oc}_X(kH))$
is a regular point.
}
\end{remark}However, the conditions on Seshadri constants are of
local nature and the results that one can deduce are strictly related to the postulation of the
chosen points. In the next section, we shall discuss one of our
results, which determines conditions on the vector bundle $\F$ and a uniform upper-bound
on the number of nodes
$\delta$ such that each point of the scheme $ {\V}_{\delta} ({\F})$ is regular.

\section{Some uniform regularity results for ${\V}_{\delta}({\E} \otimes L^{\otimes k})$}\label{S:4}

From now on, let
$X$ be a smooth projective threefold, 
$\E$ be a globally generated rank-two vector
bundle on $X$, $L$ be a very ample line bundle on $X$ and 
$k \geq 0$, $\delta >0 $ be integers.
With notation as in Section 3, we shall always take$$ {\F} = {\E} \otimes L^{\otimes k}$$and consider 
the scheme
${\V}_{\delta}({\E} \otimes L^{\otimes k})$ on $X$.
By using Theorem \ref{prop:3.fundamental} and Corollary \ref{rem:reg}, here we determine
conditions on the vector bundle $\E$ and on the integer $k$
and uniform upper-bounds on the number of nodes
$\delta$ implying that each point of 
$ {\V}_{\delta}({\E} \otimes L^{\otimes k})$ is regular. We need
before the following result.

\begin{proposition}\label{prop:7}
Let $X $ be a smooth projective threefold, $\E$ be a globally generated rank-two 
vector bundle on $X$ and L be a very ample line bundle on $X$.
Take $k > 0$ such that $L^{\otimes k}$ separates $\delta$ distinct given points
$\Sigma = \{p_1, \ldots , p_{\delta} \}$, i.e. the restriction map
\begin{equation}\label{eq:prop7}
H^0(X, L^{\otimes k}) \stackrel{\rho_k}{\to} H^0({\Oc}_{\Sigma})
\end{equation}is surjective. Thus, if $[s] \in {\V}_{\delta}({\E} \otimes L^{\otimes k})$ determines
a nodal curve $C$ in $X$ such that $Sing(C) = \Sigma$, then 
$[s] \in {\V}_{\delta}({\E} \otimes L^{\otimes k})$
is a regular point.
\end{proposition}
\begin{proof}
Since $\E$ is globally generated on $X$, the evaluation morphism
$$H^0(X, {\E}) \otimes {\Oc}_X \stackrel{ev}{\to} \E$$is surjective. This means 
that, for each $p \in X$, there exist global sections $s_1^{(p)}, \; s_2^{(p)} \in 
H^0 (X, {\E})$ such that 
$$s_1^{(p)}(p) = (1,0), \; s_2^{(p)}(p) = (0,1) \in {\Oc}_{X,p}^{\oplus 2}.$$Condition 
(\ref{eq:prop7}) means there exist global sections 
$\sigma_1, \ldots, \sigma_{\delta} \in H^0(X, L^{\otimes k})$ s. t. 
$${\sigma}_i(p_j) = \underline{0} \in {\CC}^{\delta}, \; {\rm if \; i \neq j }, \;{\rm and} \;  
\sigma_i(p_i) = (0, \ldots, \stackrel{i-th}{1}, \ldots, 0 ), \; 1 \leq i \leq 
\delta.$$Therefore, from our hypotheses, it immediately follows that 
$$H^0(X, {\E} \otimes L^{\otimes k}) \to \!\!\! \to H^0({\Oc}_{\Sigma}^{\oplus 2}) \cong 
{\CC}^{2 \delta}.$$If we take 
${\Pmc} = {\Pp}_{X}({\E} \otimes L^{\otimes k}) \stackrel{\pi}{\longrightarrow} X$ and if we 
consider, as in (\ref{eq:svolta}), $({\E} \otimes L^{\otimes k})^{\Sigma} := 
\pi_*({\Ii}_{\Sigma^1/{\Pmc}} \otimes {\Oc}_{\Pmc}(1))$, from 
diagram (\ref{eq:circledast}), we get 
\begin{displaymath}
\begin{array}{ccl}
H^0( {\E} \otimes L^{\otimes k})  & \to\!\!\! \to & H^0({\Oc}_{\Sigma}^{\oplus 2}) \cong {\CC}^{2 \delta} \\ 
\downarrow^{\mu} & & \downarrow  \\ 
H^0({\Oc}_{\Sigma}) & \stackrel{\cong}{\to} & H^0({\Oc}_{\Sigma}) \cong  {\CC}^{\delta}.\\ 
 & & \downarrow \\ 
  & & 0 
\end{array}
\end{displaymath}thus $\mu$ is surjective. By (\ref{eq:regbis}) and by the local 
description of $\alpha_X$, one can conclude.
\end{proof}

\begin{remark}\label{rem:8}
\normalfont{
With the previous result, the 
regularity condition (\ref{eq:regbis}) translates into the surjectivity 
of the restriction map $\rho_k$ in (\ref{eq:prop7}), which is a natural restriction map of 
line bundles on the threefold $X$.
}  
\end{remark}The following more general proposition gives an 
effective and uniform bound on the number $\delta = |\Sigma|$, in terms of the integer $k$, 
in order to have the surjectivity of the map $\rho_k$. 

\begin{proposition}\label{prop:9}
Let $X$ be a smooth projective $m$-fold, $L$ be a very ample line bundle and $k$ be a positive 
integer. Then, 
$L^{\otimes k}$ separates any set $\Sigma$ 
of $\delta$ distinct point of $X$ with $\delta \leq k+1$. 
In particular, 
the map $\rho_k$ in (\ref{eq:prop7}) is surjective, for each such 
$\Sigma \subset X$.
\end{proposition} 
\begin{proof}
Since $L$ is very ample on $X$, for every $p_1 \neq p_2 \in X$, there exists a section 
$s_{1,2} \in H^0(X, L)$ such that 
$$s_{1,2}(p_1) = 1 \; {\rm and} \; s_{1,2}(p_2) = 0.$$If $p_3 \in X$ is such that 
$p_3 \neq p_1, \; p_2$, there exists $s_{1,3} \in H^0 (X, L)$ such that 
$$s_{1,3}(p_1) = 1 \; {\rm and} \; s_{1,3}(p_3) = 0.$$Then 
$$\sigma : = s_{1,2} \otimes s_{1,3} \in H^0(X, L^{\otimes 2})$$is such that
$$\sigma(p_1) = 1, \; \sigma(p_2) = 0, \; \sigma(p_3) = 0.$$With analogous 
computations, it follows that $L^{\otimes 2}$ separates three points of $X$. Recursively, 
$L^{\otimes k} $ separates $k+1$ distinct points in $X$. 
\end{proof}

Finally, we have the main result of this section.
\begin{theorem}\label{thm:9bis}
Let $X$ be a smooth projective threefold, $\E$ be a globally generated rank-two vector bundle on 
$X$, $L$ be a very ample line bundle on $X$ and $k \geq 0$ and $\delta >0$ be integers. If 
\begin{equation}\label{eq:uniform}
\delta \leq k+1,
\end{equation}then ${\V}_{\delta}({\E} \otimes L^{\otimes k})$ is regular at each point.
\end{theorem}
\begin{proof}
If $k =0$, then $\delta = 1$; therefore, by the hypothesis on $\E$, it follows that 
$$H^0 (\E) \to H^0(\Oc_p^{\oplus 2}) $$is surjective, for each $p \in X$. By the definition 
of $\alpha_X$, this implies 
that ${\V}_1({\E})$ is regular at each point. 

When $k > 0$, the statement follows from Theorem \ref{prop:3.fundamental}, 
Propositions \ref{prop:7}, \ref{prop:9} and from Remark \ref{rem:8}.
\end{proof}

\begin{remark}\label{rem:10}
\normalfont{
Observe that the bound (\ref{eq:uniform}) is uniform, i.e. it does not depend 
on the postulation of nodes of the curves which are zero-loci of sections parametrized
by ${\V}_{\delta}({\E} \otimes L^{\otimes k})$. We remark that Theorem \ref{thm:9bis} generalizes
what proved by Ballico and Chiantini in \cite{BC} mainly 
because, by the characterization given in our Theorem \ref{prop:3.fundamental}, 
our approach more generally holds for families of nodal curves on smooth 
projective threefolds but also because, even in the case of $X = \Pt$, main subject 
of \cite{BC}, our result is effective and not asymptotic as Proposition 3.1 in \cite{BC}. 
Furthermore, Ballico and Chiantini
showed that in the asymptotical case, i.e. with $k >>0$, the bound $\delta \leq k+1$
is almost sharp. Indeed, they constructed an example of a non regular point
$[s] \in  {\V}_{k+4}({\Oc}_{\Pt}(k+1)\oplus {\Oc}_{\Pt}(k+4))$ whose corresponding curve
$C$ has its $(k+4)$ nodes lying on a line $L \subset \Pt$; they also showed 
that, when the points are moved so that they are in general position and 
no longer aligned, then $[s]$ is regular.
 
Thanks to Theorem \ref{thm:9bis}, the same example works not only in the asymptotic case but for 
each $k \geq 3$ proving the almost-sharpness of the bound (\ref{eq:uniform}).

} 
\end{remark}

In the next section, we also 
discuss some other examples of nodal curves on smooth projective threefolds which determine non-regular 
points of some ${\V}_{\delta}({\E} \otimes L^{\otimes k})$ (see Remarks \ref{rem:17} and \ref{rem:36}). 
  
\section{Regularity results via 
Seshadri constants and postulation of nodes}\label{S:5}

For simplicity, from now on, we focus on the 
case of $X \subset \PR$ a smooth threfold, with $L= \Oc_X(1)$ its hyperplane bundle (thus, $L^{\otimes k}$ 
will be denoted by $\Oc_X(k)$ and ${\E} \otimes L^{\otimes k}$ by ${\E}(k)$). 
As already observed in the previous section, given $X \subset \PR$ a smooth threefold, 
$\E$ a globally generated rank-two vector bundle on $X$ and $\delta >0 $, $k \geq 0$ two 
integers, Theorem \ref{thm:9bis} determines sufficient conditions in 
order that each point of the scheme ${\V}_{\delta}({\E}(k))$ is regular, for every non-negative 
integer $k$. Using a local analysis, we can determine some other regularity results which 
take into account the postulation of nodes of the curves related to the elements parametrized by 
${\V}_{\delta}({\E}(k))$.
Precisely, let $[s] \in {\V}_{\delta}({\E}(k))$, $C= V(s)$ and denote by
$\Sigma = Sing(C)$ its set of nodes. Our aim is to find some conditions on
$\Sigma$ which determine finer estimates on the 
admissible number $\delta$ of nodes in order to get the regularity of the point 
$[s] \in {\V}_{\delta}({\E}(k))$.

\begin{remark}\label{rem:13bis}
{\normalfont By Proposition \ref{prop:7} and
by Remark \ref{rem:8} a sufficient condition
for the regularity of $[s] \in {\V}_{\delta}({\E}(k))$ is to show
\begin{equation}\label{eq:van}
h^1(X, \Ii_{\Sigma/X} \otimes {\Oc}_X(k)) =0.
\end{equation}Observe that, if $\omega_X^{\vee} \otimes {\Oc}_X(k)$ is
a big and nef line bundle on $X$ then, by the Kawamata-Viehweg vanishing theorem,
(\ref{eq:van}) is exactly equivalent to the surjectivity of the map $\rho_k$ 
in \eqref{eq:prop7}, so it implies the regularity of
$[s] \in {\V}_{\delta}({\E}(k))$. In the sequel we will be concerned in finding some sufficient
conditions implying (\ref{eq:van}); we shall focus on the case when $X$ is a
Fano or a Calabi-Yau threefold and, in particular, when $X= \Pt$.
}
\end{remark}

First of all, we have to recall the following general definitions from
\cite{D}, \cite{EKL} and \cite{Ku}.

\begin{definition}\label{def:12}
Let $L$ be a nef line bundle on an $n$-dimensional projective variety $Y$. Let $p \in Y$ and 
let $b_1: Y_1 \to Y$ denote the blowing-up of $Y$ at $p$. The {\em Seshadri constant 
of $L$ at $p$}, 
$\epsilon(L,p)$, is defined as 
\begin{equation}\label{eq:sesh1}
\epsilon(L,p):= {\rm Sup} \{ \epsilon \in \RR_{\geq 0} | \; b_1^*(L) - \epsilon \; E \; {\rm 
is \; a \; nef \; {\RR}-divisor \; on} \; Y_1 \}, 
\end{equation}where $E$ denotes the $b_1$-exceptional divisor. Equivalently, 
\begin{equation}\label{eq:sesh2}
\epsilon(L,p):= {\it Inf}_{\Gamma \subset Y} \{ \frac{L \cdot \Gamma}{mult_p(\Gamma)} \},
\end{equation}where the infimum is taken over all reduced and irreducible curves 
$\Gamma \subset Y$ passing through $p$.

\noindent
More generally, if $\delta$ is an integer greater than $1$ and if $p_1, \; \ldots,\; p_{\delta} 
\in Y$ are $\delta$ distinct points then, denoting by $b_{\delta} : Y_{\delta} \to Y$ the 
blowing-up of $Y$ along the given points, the {\em multiple point Seshadri constant at 
$p_1, \; \ldots,\; p_{\delta}$} is defined as 
\begin{equation}\label{eq:msesh1}
\epsilon(L, p_1, \; \ldots,\; p_{\delta}):= 
{\rm Sup} \{ \epsilon \in \RR_{\geq 0} | \; b_{\delta}^*(L) - \epsilon \; \Sigma_{i=1}^{\delta}
E_i \; {\rm 
is \; a \; nef \; {\RR}-divisor \; on} \; Y_{\delta} \}, 
\end{equation}where $\Sigma_{i=1}^{\delta}E_i$ is the 
$b_{\delta}$-exceptional divisor. As before, one also has 
\begin{equation}\label{eq:msesh2}
\epsilon(L, p_1, \; \ldots,\; p_{\delta}):= {\it Inf}_{\Gamma \subset Y} 
\{ \frac{L \cdot \Gamma}{\Sigma_{i=1}^{\delta}mult_{p_i}(\Gamma)} \},
\end{equation}where the infimum is taken over all integral curves 
$\Gamma \subset Y$ s.t. $\Gamma \cap \{ p_1, \; \ldots, \; p_{\delta} \} \neq \emptyset$.
\end{definition}

\begin{definition}\label{def:13}
Let $Y$ be a projective variety of dimension $n$ and $\delta \geq 2$ be a positive 
integer. Let $Y^{(\delta)}$ denote the $\delta$-Cartesian product of $Y$ minus the diagonals. 
If $(p_1, \; \ldots,\; p_{\delta}) \in Y^{(\delta)}$, the points $p_1, \; \ldots,\; p_{\delta}$ 
are called {\em general points of} $Y$ if $(p_1, \; \ldots,\; p_{\delta})$ is 
outside a Zarisky closed subset of 
$Y^{(\delta)}$ and {\em very general points of} $Y$ if $(p_1, \; \ldots,\; p_{\delta})$ is 
outside the union of countably many proper subvarieties of $Y^{(\delta)}$.
\end{definition}

\noindent
Before stating our next result, we recall that from our assumptions in \S 2 (see Definition 
\ref{def:0} and above) the integer $\delta$ is always assumed to be 
$\delta \leq min \{h^0({\E}(k)) -1, p_a(C) \}$.

\begin{theorem}\label{thm:14}
Let $X \subset \PR$ be a smooth threefold such that $deg(X) = d$, $\omega_X \cong \Oc_X(-m)$ 
for some integers $d>0$ and $m \geq 0$. Let $\E$ be a globally generated rank-two vector 
bundle on $X$ and let $k \geq 0$ and $\delta >0$ 
be integers. Let $[s] \in {\V}_{\delta}({\E}(k))$, $C= V(s)$ and let
$\Sigma$ denote its set of nodes.
Assume either

\noindent
(i) $k+m>3$, when $\delta=1$, or

\noindent
(ii) $\Sigma$ is a set of $\delta \geq 2$ very general points on $X$ and

\begin{itemize}
\item[a)]$k+m > \frac{6}{\sqrt[3]{d}}$, when $d < 8$ and $\delta \leq 5$;
\item[b)] $k+m> max \{\frac{36}{deg(C)}, \frac{18}{\sqrt[3]{25d}} \}$, when $d<8$ and
$6 \leq \delta < min \{h^0({\E}(k)) , \frac{1}{6}(k+m) deg(C), \; \delta_0^{(k)}\}$, where
$\delta_0^{(k)}$ is a root of the polynomial $F_{k,m,d}(\delta) := 27 \delta^3 -
(k+m)^3 d (\delta-1)^2$ such that $F_{k,m,d}(\delta) < 0$ on the connected interval
$[6,  \delta_0^{(k)})$;
\item[c)] $k+m >3$, when $d\geq 8$ and $\delta \leq d-2$;
\item[d)] $k+m > max \{\frac{6(d-2)}{deg(C)}, \frac{3(d-2)}{\sqrt[3]{d(d-2)^2}} \}$ when
$d \geq 8$ and $d-2 \leq \delta < min \{h^0({\E}(k)) , \frac{1}{6}(k+m) deg(C), \; \delta_0^{(k)}\}$, where
$\delta_0^{(k)}$ is a root of the polynomial $F_{k,m,d}(\delta) := 27 \delta^3 -
(k+m)^3 d (\delta-1)^2$ such that $F_{k,m,d}(\delta) < 0$ on the connected interval
$[d-2,  \delta_0^{(k)})$.
\end{itemize}
Then, in each case, $[s]$ is a regular point of ${\V}_{\delta}({\E}(k))$.
\end{theorem}
\begin{proof}
Let $b_{\delta} : Y_{\delta} \to X$ be the blowing-up of $X$ along $\Sigma$. From our
assumptions on $X$ and from Leray's isomorphism, it follows that
$$H^1( X, {\Ii}_{\Sigma/X}(k)) = H^1( X, {\Ii}_{\Sigma/X}(k+m) \otimes \omega_X) \cong
H^1(Y_{\delta}, \omega_{Y_{\delta}} \otimes \Oc_{Y_{\delta}}((k+m) b_{\delta}^*(H) - 3 B)),$$
where $B= \sum_{i=1}^{\delta} E_i$ is the $b_{\delta}$-exceptional divisor.
Therefore, if $(k+m) b_{\delta}^*(H) - 3 B$ is a big and nef divisor, by the Kawamata-Viehweg
vanishing theorem, $H^1( X, {\Ii}_{\Sigma/X}(k)) = (0)$, which implies the regularity of
$[s] \in {\V}_{\delta}({\E}(k))$ (see Remark \ref{rem:13bis}).

Let $\overline{\epsilon} =
\epsilon ({\Oc}_X(1), \Sigma)$ denote the multiple point Seshadri constant of the very ample line
bundle ${\Oc}_X(1)$ at $\Sigma$; then $$(k+m) b_{\delta}^*(H) - 3 B =
\frac{3}{\overline{\epsilon}}( b_{\delta}^*(H) -  \overline{\epsilon} B) +
(k+m-\frac{3}{\overline{\epsilon}}) b_{\delta}^*(H).$$Observe that the first summand in the right
hand side is nef, by definition of $\overline{\epsilon}$, whereas the second is big and nef
as soon as $\overline{\epsilon} > \frac{3}{k+m}$.

We want to show that our hypotheses imply that the Seshadri constant $\overline{\epsilon}$ is
always greater than $\frac{3}{k+m}$; so the statement will be proved.

\noindent
(i) If $\delta = 1$, then $\overline{\epsilon}= \epsilon ({\Oc}_X(1),p) \geq 1$, for each
$p \in X$, since ${\Oc}_X(1)$ is very ample. Therefore, since $k+m >3$ implies $\frac{3}{k+m}<1$, we
have $h^1(X, {\Ii}_{\{p\}/X}(k)) =0$, for each $p \in X$.

\noindent
(ii) For $\delta \geq 2$, we can consider Theorem 1.1 in \cite{Ku}. For $L$ a big and nef
line bundle on $X$, the author denotes by $\epsilon (L; \delta) $ the Seshadri constant of $L$
at very general $\delta$ points of $X$, whereas, by $\epsilon (L; 1)$ the Seshadri constant of $L $ at a
very general point of $X$. In the threefold case with $L={\Oc}_X(1)$,
K$\ddot{u}$chle's result
gives
\begin{equation}\label{eq:nettuno}
\epsilon ({\Oc}_X(1); \delta) \geq M:= {\rm min} \{ \epsilon ({\Oc}_X(1); 1), \;
\frac{\sqrt[3]{d}}{2}, \;
\frac{\sqrt[3]{d(\delta-1)^2}}{\delta}\},
\end{equation}
where $\epsilon ({\Oc}_X(1); 1) \geq 1$, since ${\Oc}_X(1)$ is very
ample. By assumption, $\Sigma$ is a set of very general points on $X$, thus
$\overline{\epsilon} = \epsilon({\Oc}_X(1); \Sigma)$ coincides with $\epsilon({\Oc}_X(1);\delta)$.
Therefore, to prove that $\overline{\epsilon} > \frac{3}{k+m}$ we reduce to showing that
our numerical hypotheses imply
\begin{equation}\label{eq:nettuno1}
\frac{3}{k+m} < M.
\end{equation}Observe that, when $d \geq 8$ and $\delta \leq d-3$, we have $M \geq 1$, since all
the real numbers in the brackets in (\ref{eq:nettuno}) are greater than or equal to $1$. Since $k+m >3$,
then (\ref{eq:nettuno1}) trivially holds.

\noindent
In the other cases, we find that:

\begin{itemize}
\item $\epsilon({\Oc}_X(1) ; 1)$ is always greater than or equal to $1$, since ${\Oc}_X(1)$ is
very ample;
\item $\frac{\sqrt[3]{d}}{2} <1$ iff $d<8$;
\item $\frac{\sqrt[3]{d(\delta-1)^2}}{\delta} <1$ if $\delta \geq d-2$, when $d \geq 3$,
or if $\delta \geq 2$, when $1 \leq d \leq 2$;
\item $\frac{\sqrt[3]{d(\delta-1)^2}}{\delta} <\frac{\sqrt[3]{d}}{2} $ iff $\delta \geq 6$.
\end{itemize}Therefore, considering all the above inequalities, we find that

\[ M = \left\{ \begin{array}{lcl}
\frac{\sqrt[3]{d}}{2}& {\rm if} &  d < 8 \; {\rm and} \; 2 \leq \delta \leq 5,\\
\frac{\sqrt[3]{d(\delta-1)^2}}{\delta} & {\rm if} & {\rm either} \; d \geq 8 \;
                                                     {\rm and} \; \delta \geq d-2 \\
                           &         & {\rm or} \; d < 8 \; {\rm and} \; \delta \geq 6.
                           \end{array}
\right. \]In all these cases we have $M <1$.

\noindent
When $M=\frac{\sqrt[3]{d}}{2}$, (\ref{eq:nettuno1}) holds as soon as
$k+m > \frac{6}{{\sqrt[3]{d}}}$. On the other hand,
when $M= \frac{\sqrt[3]{d(\delta-1)^2}}{\delta}$, we want
\begin{equation}\label{eq:*}
\frac{3}{k+m} < \frac{\sqrt[3]{d(\delta-1)^2}}{\delta}.
\end{equation}Since this case occurs when
$d \geq 8$, $\delta \geq d-2$ and when $d<8$, $\delta \geq 6$, we impose
\begin{equation}\label{eq:**1}
\frac{3(d-2)}{\sqrt[3]{d(d-3)^2}} < k+m, \; {\rm when} \; d \geq 8,
\end{equation}and
\begin{equation}\label{eq:**2}
\frac{18}{\sqrt[3]{25d}} < k+m, \; {\rm when} \; d < 8.
\end{equation}Observe
that (\ref{eq:*}) is equivalent to asking that the polynomials
$$F_{k,m,d}(\delta) := 27 \delta^3 - d (k+m)^3 (\delta -1)^2$$satisfy the inequalities
$F_{k,m,d}(\delta) <0$. By (\ref{eq:**1}) and (\ref{eq:**2}), we have that
$$F_{k,m,d}(d-2) < 0, \: {\rm when} \; d \geq 8,$$
$$F_{k,m,d}(6) < 0, \: {\rm when} \; d < 8.$$Therefore
each cubic polynomial $F_{k,m,d}(\delta)$ has (at least) one root
which is greater than $d-2$, when $d \geq 8$, and greater than $6$, when $d<8$,
respectively. Denote by $\delta_0^{(k)}$ the root of $F_{k,m,d}(\delta)$ s.t.
$$F_{k,m,d}(\delta) <0 , \; \forall \; \delta \in [d-2,\delta_0^{(k)}), \; {\rm when} \; d \geq 8,$$
$$F_{k,m,d}(\delta) <0 , \; \forall \; \delta \in [6,\delta_0^{(k)}), \; {\rm when} \; d < 8,$$
respectively. Thus, in such ranges of values for $\delta$, (\ref{eq:nettuno1}) automatically
holds.

On the other hand, since $[s] \in {\V}_{\delta}({\E})$,
then, $C \subset X$ is an
irreducible curve having nodes at $\Sigma$; thus, by definition of multiple point
Seshadri constant - see
(\ref{eq:msesh2}) - we have $\frac{deg(C)}{mult_{\Sigma}(C)} > \frac{3}{k+m}$, i.e.
$\delta < \frac{1}{6} (k+m) deg(C)$. Therefore,
when $\delta \geq d-2$, we have $(k+m) > \frac{6(d-2)}{deg(C)}$, whereas
$\delta \geq 6$ gives $(k+m) > \frac{36}{deg(C)}$.
\end{proof}

\noindent
When, in particular, $X = \Pt$ we can simplify the previous result.

\begin{coroll}\label{cor:14bis}
Let $\E$ be a globally generated rank-two vector bundle on $\Pt$. Denote by $c_i$ the
$i^{th}$-Chern class of $\E$. 
Let $k$ and $\delta$ be integers such that $k \geq 0$ and $\delta >0$. 
Let $[s] \in {\V}_{\delta}({\E}(k))$ and let 
$\Sigma$ denote the set of nodes of the curve $C \subset \Pt$ corresponding to $s$. Assume that 

\noindent
(i) $k\geq 0$, when $\delta=1$;

\noindent
(ii) $k \geq 3$, when:  

\begin{itemize}
\item[a)] $\delta = 2$,
\item[b)] $3 \leq \delta \leq 5$ and $\Sigma$ is a set of very general points in $\Pt$, 
\item[c)] $3 c_1 + c_2 + 4 >0$, $\Sigma$ is a set of very general points in $\Pt$ and 
$6 \leq \delta < min \{h^0({\E}(k)) , \frac{1}{6}(k+4)(k^2 + c_1 k + c_2) , \; \delta_0^{(k)}\}$, where
$\delta_0^{(k)}$ is a positive root of the polynomial
$F_{k,}(\delta) := 27 \delta^3 - (k+4)^3 (\delta-1)^2$ such that $F_{k}(\delta) < 0$ on
the connected interval
$[6,  \delta_0^{(k)})$.
\end{itemize}
Then, in each case, $[s]$ is a regular point of ${\V}_{\delta}({\E}(k))$.
\end{coroll}
\begin{proof}
As in the proof of Theorem \ref{thm:14}, it suffices to show that the multiple
point Seshadri constant $\overline{\epsilon} := \epsilon({\Oc}_X(1); \; \Sigma) > \frac{3}{k+4}$.

\noindent
If $\delta =1$, $\overline{\epsilon}= \epsilon ({\Oc}_X(1); \; p) =1 $ for each $p \in \Pt$,
since there exist lines in $\Pt$. Thus, for each $k \geq 0$ we have
$\overline{\epsilon} > \frac{3}{k+4}$;

\noindent
If $\delta =2$, then $\overline{\epsilon}= \epsilon ({\Oc}_X(1); \; p_1, p_2) \geq
\frac{1}{2}$ for
$(p_1, p_2) \in (\Pt)^{(2)}$ (as in Definition \ref{def:13}), since for all $p_1 \neq p_2 $ there exists
the line
$L_{p_1,p_2} = < p_1, p_2>$. Therefore, if $k \geq 3$,
$\overline{\epsilon} > \frac{3}{k+4}$;

\noindent
For $\delta \geq 3$ we can use the same procedure of Theorem \ref{thm:14} observing that
$\epsilon({\Oc}_X(1); \Sigma) = \epsilon({\Oc}_X(1); \delta)$, by assumption on $\Sigma$, and
that $ \epsilon ({\Oc}_X(1) , \delta) \geq M:= {\rm min} \{ 1, \; \frac{1}{2},\;
\frac{\sqrt[3]{(\delta-1)^2}}{\delta} \}.$
\end{proof}

\begin{remark}\label{rem:15}
\normalfont{
Observe that the polynomials $F_k(\delta)$ in Corollary \ref{cor:14bis}
asymptotically give the upper-bounds
$ \delta< \frac{(k+4)^3}{27}$ (equivalently $k > 3 \sqrt[3]{\delta} -4$ ). Therefore,
we have a cubic polynomial in the indeterminate $k$ which bounds the admissible
number of nodes of $C$. The same occurs with the inequalities $\delta < h^0(X, {\E}(k))$ and
$\delta< \frac{1}{6}(k+4)(k^2+c_1k + c_2)$. Similar situation for the polynomials
$F_{k,m,d}(\delta)$ in Theorem \ref{thm:14}. Therefore, we have
cubic upper-bounds on $k$ for
$\delta$ to get regularity results for the point
$[s] \in {\V}_{\delta}({\E}(k))$. This depends on the fact that the computations are
related to Seshadri constants of very ample line bundles at
very general points on a $3$-dimensional variety.
}
\end{remark}

\begin{remark}\label{rem:15bis}
\normalfont{
At this point, on the one hand we have 
Theorem \ref{thm:9bis} which gives uniform bounds on the admissible number $\delta$ 
of nodes in order that each point of the Severi variety ${\V}_{\delta}({\E}(k))$ on a threefold
$X$ is regular; morover, these uniform upper-bounds
only depend on the number of nodes and not on their configurations in $X$. On the other hand,
Theorem \ref{thm:14} determines sufficient conditions for the
regularity of a point $[s] \in {\V}_{\delta}({\E}(k))$, having assumed that
the nodes of $C=V(s)$ are in very general position on $X$. Therefore, there are
intermidiate cases which are very interesting to study. Precisely, if 
$[s] \in {\V}_{\delta}({\E}(k))$, we want to
find conditions for its regularity 
assuming that a (not necessarily proper) subset of the nodes of $C$
lies on a proper subscheme of $X$.
}
\end{remark}

Given $X \subset \PR$ a smooth threefold,
$\E$ a globally generated rank-two vector bundle on $X$
and $k$ and $\delta$ positive integers, consider $[s]\in {\V}_{\delta}({\E}(k))$.
From now on, $C$ will denote the irreducible curve determined by $s$, whose set of nodes is
$\Sigma= Sing(C)$, as well as $\Sigma_0 \subseteq \Sigma$ will denote a
(not necessarily proper) subset of nodes of $C$.

\begin{proposition}\label{prop:16}
Let $X \subset \PR$ be a smooth threefold and let $H$ denote
its hyperplane section. Let $S_a \subset X$ be an irreducible
divisor such that $S_a \sim aH$ on $X$, for some positive integer $a$.
Let $[s] \in {\V}_{\delta}({\E}(k))$, $C = V(s)$
and assume that $\Sigma_0 \subset \Sigma=Sing(C)$ lies
on $S_a \setminus Sing(S_a)$. Assume also that $$h^1(X, {\Oc}_X(k-a)) = h^1(X, {\Oc}_X(k)) =
h^2(X, {\Oc}_X(k-a)) =0$$(e.g, when $X$ is arithmetically Cohen-Macaulay). Then,
if $\Sigma_0$ does not impose independent conditions to the complete linear system
in $|{\Oc}_{S_a}(k)|$ on $S_a$, $[s]$ cannot be a regular point for
${\V}_{\delta}({\E}(k))$.

\end{proposition}

\begin{proof}
Since $S_a \sim aH$ on $X$, we have the ideal sequence
\begin{equation}\label{eq:ideal1}
0 \to {\Oc}_X(k-a) \to {\Ii}_{\Sigma_0/X}(k) \to {\Ii}_{\Sigma_0/S_a}(k) \to 0.
\end{equation}Therefore, by the hypotheses on $X$,
$H^1(X, {\Ii}_{\Sigma_0/X}(k)) \cong H^1(S_a, {\Ii}_{\Sigma_0/S_a}(k))$ which
implies the statement.
\end{proof}In particular,

\begin{coroll}\label{cor:16bis}
Let $S_a \subset \Pt$ be an irreducible (not necessarily smooth)
surface of degree $a$. Let $[s] \in {\V}_{\delta}({\E}(k))$, $C=V(s)$ and assume that $\Sigma_0 \subseteq
\Sigma =Sing(C)$ is such that $\Sigma_0 \subset S_a \setminus Sing (S_a)$. Then, if
$\Sigma_0$ does not impose
independent conditions to the complete linear system $|{\Oc}_{S_a}(k)|$
on $S_a$, $[s]$ is not a regular point for ${\V}_{\delta}({\E}(k))$.
\end{coroll}

\begin{remark}\label{rem:17}
\normalfont{
Observe that, with Proposition \ref{prop:16} and Corollary \ref{cor:16bis}, one can
easily construct many examples of non-regular points $[s] \in {\V}_{\delta}({\E}(k))$,
corresponding to nodal curves on a smooth
projective threefold $X$, by translating the problem to linear systems on surfaces $S \subset X$
not separating a given set of smooth points in $S$.

On the other hand, one can
also find some conditions on $\delta_0 = |\Sigma_0|$ ensuring that, if
$[s] \in {\V}_{\delta}({\E}(k))$ is not regular, the failure of the regularity property
depends on the behaviour of the nodes in $\Sigma \setminus \Sigma_0$. Indeed, since $\Sigma_0 \subseteq \Sigma$,
we have
\begin{equation}\label{eq:(*)}
0 \to {\Ii}_{\Sigma/X}(k) \to {\Ii}_{\Sigma_0/X}(k) \to
{\Oc}_{\Sigma \setminus \Sigma_0}(k) \to 0.
\end{equation}Taking $X$ as in Proposition \ref{prop:16} and assuming that
$h^1(S_a, {\Oc}_{S_a}(k))=0$ (e.g, for $X=\Pt$), then if we have
some conditions implying that $|{\Oc}_{S_a}(k)|$ separates $\Sigma_0$ on $S_a$, by
(\ref{eq:(*)}) we have $$H^0({\Oc}_{\Sigma \setminus \Sigma_0}(k)) \to \!\!\to
H^1({\Ii}_{\Sigma/X}(k)).$$Therefore, a possibly non-zero element in
$H^1(X, {\Ii}_{\Sigma/X}(k))$ is induced by an element in
$H^0({\Oc}_{\Sigma \setminus \Sigma_0}(k))$.
}
\end{remark}To get some effective results, we can use several approaches. First of all, 
we want to consider the case when $S_a \subset X$ is a smooth irreducible divisor in $X$, which 
is linearly equivalent to $aH$ on $X$. As in Theorem \ref{thm:14}, our aim is 
to show that if $\Sigma_0 \subset S_a$ is a set of very general points 
on $S_a$, then we get a quadratic upper-bound on $k$ for the admissible number 
$\delta_0 = |\Sigma_0|$ in order that the complete linear system $|{\Oc}_{S_a}(k)|$ separates 
$\Sigma_0$ on $S_a$.

\begin{proposition}\label{prop:18}
Let $X\subset \PR$ be a smooth threefold of degree $d$ and let $H$ 
denote its hyperplane section. Assume that $\omega_X \cong {\Oc}_X(-m)$, 
for some integer $m \geq 0$. Let $S_a$ be a smooth irreducible divisor such that 
$S_a \sim aH$ on $X$, for some positive integer $a$. Denote by $H_{S_a}$ the hyperplane 
section of $S_a \subset \PR$. Let $\Sigma_0\subset S_a$ be a set of $\delta_0$ 
distinct points on $S_a$. Given $k$ a non-negative integer, assume that:

\noindent
(i) $k+m > a+2$, when either
\begin{itemize}
\item[a)] $\delta_0=1$, or
\item[b)] $ad \geq 4$ and $2 \leq \delta_0 < \frac{ad + \sqrt{ad(ad-4)}}{2}$, 
\end{itemize}

\noindent
(ii) $k+m > a + \frac{4}{\sqrt{ad}}$, when either
\begin{itemize}
\item[a)] $ad < 4$ and $2 \leq \delta_0 < \frac{(k+m-a)^2 + 
\sqrt{(k+m-a)^2ad((k+m-a)^2ad-16)}}{8}$, or
\item[b)] $ad \geq 4$ and $\frac{ad + \sqrt{ad(ad-4)}}{2} \leq \delta_0 
< \frac{(k+m-a)^2 + \sqrt{(k+m-a)^2ad((k+m-a)^2ad-16)}}{8}$. 
\end{itemize}Then, 
\begin{equation}\label{eq:(**)}
h^1(S_a, {\Ii}_{\Sigma_0/S_a}(k))=0.
\end{equation}
\end{proposition}
\begin{proof}
Let $b_{\delta_0} : \tilde{S}_a \to S_a$ be the blowing-up of $S_a$ along $\Sigma$. Then 
\begin{displaymath}
\begin{array}{cl}
H^1(S_a, {\Ii}_{\Sigma_0/S_a}(k)) = & H^1(S_a, {\Ii}_{\Sigma_0/S_a}(k+m-a) \otimes 
\omega_{S_a})\cong \\
 & H^1(\tilde{S}_a, \omega_{\tilde{S}_a} \otimes {\Oc}_{\tilde{S}_a}(
(k+m-a) b_{\delta_0}^*(H_{S_a}) - 2B))\;\;\;(*),
\end{array}
\end{displaymath}where $B=\sum_{j=1}^{\delta_0} E_j$ is the 
$b_{\delta_0}$-exceptional divisor. Let $\tilde{\epsilon} = \epsilon({\Oc}_{S_a}(1), \Sigma_0)$ 
be the multiple point Seshadri constant of ${\Oc}_{S_a}(1)$ at $\Sigma_0$, then 
$$(k+m-a) b_{\delta_0}^*(H_{S_a}) - 2B = 
\frac{2}{\tilde{\epsilon}}(b_{\delta_0}^*(H_{S_a}) - \tilde{\epsilon}B) +
(k+m-a - \frac{2}{\tilde{\epsilon}}) b_{\delta_0}^*(H_{S_a}).$$Therefore, if 
$\tilde{\epsilon} > \frac{2}{k+m-a}$ (with $a \neq k+m$) by the Kawamata-Viehweg 
vanishing theorem we get the desired vanishing in $(*)$. 

\noindent
If $\delta_0=1$, then $\epsilon ({\Oc}_{S_a}(1); p) \geq 1$ for each $p \in S_a$, since 
${\Oc}_{S_a}(1)$ is very ample. Therefore, if $k+m > a+2$, the vanishing in $(*)$ holds. 

\noindent
If $\delta_0 \geq 2$, by Theorem 1.1 in \cite{Ku}, if 
$\epsilon ({\Oc}_{S_a}(1); \delta_0)$ 
denotes the multiple point Seshadri constant of ${\Oc}_{S_a}(1)$ at $\delta_0$ 
very general points of 
$S_a$, then 
$$\epsilon ({\Oc}_{S_a}(1); \delta_0) \geq M:= {\rm min} \{\epsilon ({\Oc}_{S_a}(1); 1), 
\frac{\sqrt{ad}}{2}, \frac{\sqrt{ad(\delta_0-1)}}{\delta_0} \}.$$By straightforward 
computations, if 
$\Sigma_0 \subset S_a$ is a set of $\delta_0$ very general points on $S_a$ and 
if our numerical hypotheses hold, then the vanishing in $(*)$ holds.
\end{proof}

As a consequence of Proposition \ref{prop:16}, Remark \ref{rem:17} and Proposition 
\ref{prop:18}, we have the following:

\begin{theorem}\label{thm:19}
Let $X\subset \PR$ be a smooth threefold of degree $d$ and let ${\Oc}_X(1)$ denote 
its hyperplane bundle. Assume that 
$\omega_X \cong {\Oc}_X(-m)$, for some $m \geq 0$. Let $S_a$ be a smooth 
irreducible divisor on $X$, such that $S_a \in |{\Oc}_X(a)|$ on $X$, for some $a >0$. Let 
$[s] \in {\V}_{\delta}({\E}(k))$, where $\E$ is a globally generated rank-two 
vector bundle on $X$ and $\delta$ and $k$ two positive integers. If $C = V(s)$ and if 
$\Sigma = Sing(C)$, assume that $\Sigma_0 \subseteq \Sigma$ lies on $S_a$. Assume further that 
$$H^1(X, {\Oc}_X(k-a)) = H^1(X, {\Oc}_X(k)) = H^2(X, {\Oc}_X(k-a))= (0)$$and that the 
numerical hypotheses in Proposition \ref{prop:18} hold. 
Then, if $[s] \in {\V}_{\delta}({\E}(k))$ is not a regular point, the failure of the 
regularity property depends on the points in $\Sigma \setminus \Sigma_0$. 

\noindent
In particular, if $\Sigma = \Sigma_0$, $[s] \in {\V}_{\delta}({\E}(k))$ is a regular point. 
\end{theorem}
\begin{proof}From our assumptions, by Proposition \ref{prop:16} we get 
that $$H^1(X, {\Ii}_{\Sigma_0/X}(k)) \cong H^1(X, {\Ii}_{\Sigma_0/S_a}(k)).$$Now, from 
Proposition \ref{prop:18} and from the exact sequence (\ref{eq:(*)}), it 
follows that $$H^0({\Oc}_{\Sigma \setminus \Sigma_0}(k) ) \to \!\! \to H^1(X, 
{\Ii}_{\Sigma/X}(k)),$$i.e. if there exists a non-zero obstruction, it is induced by an element 
in $H^0({\Oc}_{\Sigma \setminus \Sigma_0}(k))$. 
\end{proof}When $X=\Pt$, the above results reduce to:

\begin{coroll}\label{cor:20}
Let $\E$ be a globally generated rank-two vector bundle on
$\Pt$ and let $k$ and $\delta$ be positive integers. Let $[s] \in {\V}_{\delta}({\E}(k))$
and let $C= V(s)$ with $\Sigma= Sing(C)$. Assume that $\Sigma_0 \subseteq \Sigma$ lies on
a smooth surface $S_a$ of degree $a$ and let $\delta_0$ be the cardinality of $\Sigma_0$. 
Assume that the points in $\Sigma_0$ are in very general position on $S_a$ and that the 
following conditions hold:

\noindent
(i) $k > a - 2$, when either
\begin{itemize}
\item[a)] $ a \geq 1$ and $\delta_0=1$, or
\item[b)] $a \geq 4$ and $2 \leq \delta_0 < \frac{a + \sqrt{a(a-4)}}{2}$, 
\end{itemize}

\noindent
(ii) $k> a - 4 + \frac{4}{\sqrt{a}}$, when either
\begin{itemize}
\item[a)] $a < 4$ and $2 \leq \delta_0 < \frac{(k+4-a)^2 a + 
\sqrt{(k+4-a)^2a((k+4-a)^2a-16)}}{8}$, or
\item[b)] $a \geq 4$ and $\frac{a + \sqrt{a(a-4)}}{2} \leq \delta_0 
< \frac{(k+4-a)^2 a + 
\sqrt{(k+4-a)^2a((k+4-a)^2a-16)}}{8}$. 
\end{itemize}Then, $H^1({\Pt}, {\Ii}_{\Sigma_0/{\Pt}}(k)) = (0)$. In particular, if 
$[s] \in {\V}_{\delta}({\E}(k))$ fails to be a regular 
point, the failure of the regularity property depends on the 
nodes in $\Sigma \setminus \Sigma_0$. In particular, if 
$\Sigma = \Sigma_0$, $[s] \in {\V}_{\delta}({\E}(k))$ is a regular point.
\end{coroll}

Since we are interested in very general points, we can generalize 
the previous approach by assuming that $S_a=S$ is not necessarily smooth. Indeed in general, 
if $\mu: \tilde{S} \to S$ denotes a resolution of singularities for 
$S$, given $L$ a Weil divisor on $S$, we have
\begin{equation}\label{eq:nettuno2} 
\epsilon (L; \mu(p_1), \cdots, \mu(p_{\delta})) = \epsilon (\mu^*(L); p_1, \ldots, 
p_{\delta}),
\end{equation}since around the $p_i$'s $\mu$ is an isomorphism. Therefore, we can more generally 
consider $S$ to be a normal surface with sufficiently mild singularities.

For simplicity, we shall discuss the case of $S \subset \Pt$ of degree $a$; the case
$S \subset X \subset \PR$, where $X$ a smooth threefold, is a straightforward generalization.

\begin{theorem}\label{thm:21}
Let $S \subset \Pt$ be an irreducible surface of degree
$a$ having at worst log-terminal singularities and denote by $H$ the hyperplane section of $S$.
Suppose that $K_S$ is a $\QQ$-Cartier divisor of index $r$, such that $r K_S \equiv
\alpha H$, for some $\alpha \in \NN$. Let $\E$ be a globally generated rank-two vector bundle
on $\Pt$, $k$ a positive integer and let $[s] \in {\V}_{\delta}({\E}(k))$. Let $C= V(s)$ and let
$\Sigma = Sing(C)$. Assume that $\Sigma_0 \subseteq \Sigma$ is a set of very
general points in $S$. Assume further that

\noindent
(i) $k > 2 + \frac{\alpha}{r}$, when either
\begin{itemize}
\item[a)] $ a \geq 1$ and $\delta_0=1$, or
\item[b)] $a > 4$ and $2 \leq \delta_0 < \frac{a + \sqrt{a(a-4)}}{2}$, 
\end{itemize}

\noindent
(ii) $k > \frac{4}{\sqrt{a}} + \frac{\alpha}{r}$, when $a \leq 4$ and $\delta_0=2$,

\noindent
(iii) $k> \frac{4r}{\sqrt{a}} + \frac{\alpha}{r}$, when either
\begin{itemize}
\item[a)] $a \leq 4$ and $3 \leq \delta_0 < \frac{a(rk-\alpha)^2  + 
\sqrt{a(rk-\alpha)^2 ((a(rk-\alpha)^2 -16r^2)}}{8r^2}$, or
\item[b)] $a > 4$ and $\frac{a + \sqrt{a(a-4)}}{2} \leq \delta_0 
< \frac{a(rk-\alpha)^2  + 
\sqrt{a(rk-\alpha)^2 ((a(rk-\alpha)^2 -16r^2)}}{8r^2}$. 
\end{itemize}Then, if 
$[s] \in {\V}_{\delta}({\E}(k))$ is not a regular point, the failure of the regularity 
property depends on the nodes in $\Sigma \setminus \Sigma_0$. In particular, if 
$\Sigma = \Sigma_0$, $[s]\in {\V}_{\delta}({\E}(k))$ is a regular point.
\end{theorem} 
\begin{proof}
Let $\mu : Y \to S$ be a log-resolution of the pair $(S,0)$ (see $\S$ 1). Then 
$$K_Y + \Delta \equiv \mu^*(K_S) +P,$$
where $\Delta$ is a boundary divisor and $P$ is an integral, effective and $\mu$-exceptional 
divisor on $Y$. Since $\Sigma_0 \subset S$ is a set of very general point 
on $S$, then $\mu^*(\Sigma_0) = \Sigma_0' \cong \Sigma_0$ is a set of very general points 
on $Y$. 
Let $b_{\delta_0}: \tilde{Y} \to Y$ be the blowing-up of $Y$ along $\Sigma_0'$ and denote by 
$F_{\delta_0} : \tilde{Y} \to S$ the composition $F_{\delta_0} = b_{\delta_0} \circ \mu$. 
Since, by hypothesis, $rK_S \equiv \alpha H$ then 
\begin{displaymath}
\begin{array}{rll}
h^1(S, {\Ii}_{\Sigma_0/S}(k)) & = & h^1(S, {\Ii}_{\Sigma_0/S}((k-\frac{\alpha}{r})H + K_S) ) \\
 & = &h^1(Y, {\Ii}_{\Sigma_0'/Y} \otimes {\Oc}_Y(((k-\frac{\alpha}{r})\mu^*(H) + K_Y + 
\Delta) )) \\
 & = & h^1(\tilde{Y}, {\Oc}_{\tilde{Y}}(K_{\tilde{Y}} + b_{\delta_0}^*(\Delta)+ 
(k-\frac{\alpha}{r})F_{\delta_0}^*(H) -2B)),
\end{array}
\end{displaymath}where $B= \sum_{i=1}^{\delta_0} E_i$ is the $b_{\delta_0}$-exceptional 
divisor. Since, by (\ref{eq:nettuno2}), 
$\overline{\epsilon} = \epsilon (\mu^*({\Oc}_S(1)), \Sigma_0') = \epsilon({\Oc}_S(1), 
\Sigma_0)$, then 
$$ (k-\frac{\alpha}{r})F_{\delta_0}^*(H) -2B = (k-\frac{\alpha}{r} - 
\frac{2}{\overline{\epsilon}})F_{\delta_0}^*(H) + \frac{2}{\overline{\epsilon}} 
(F_{\delta_0}^*(H) - \overline{\epsilon}B)$$is big and nef if 
$k-\frac{\alpha}{r} > \frac{2}{\overline{\epsilon}}$, i.e. if ${\overline{\epsilon}} 
> \frac{2r}{rk - \alpha}$. At this point, we can apply the same 
computations as in Proposition \ref{prop:18} and in Theorem \ref{thm:19}.
\end{proof}

\begin{remark}\label{rem:22}
\normalfont{
As observed in Remark \ref{rem:15} for sets of points in very general position on a threefold, 
from Corollary \ref{cor:20} and Theorem \ref{thm:21} we see that, 
when $\Sigma_0 = \Sigma$ is assumed to be a set of very general points lying on a 
smooth surface in $\Pt$ or on a normal surface with at worst log-terminal singularities, 
there are some upper-bounds on the number of admissible points in $\Sigma$ such that an element 
$[s] \in {\V}_{\delta}({\E}(k))$, whose zero-locus $C$ has nodes at $\Sigma$, is a regular point. 
Such upper-bounds are quadratic polynomials in $k$; this reflects 
the fact that $\Sigma$ is assumed to be a set of very 
general points lying on a $2$-dimensional subscheme. Indeed, all the computations 
are related to multiple point Seshadri constants of very ample line bundles on such schemes.
}
\end{remark}

Looking back at Proposition \ref{prop:18}, at Theorem \ref{thm:21} and at
Corollary \ref{cor:20}, we want to study how the upper-bounds on $\delta_0 = |\Sigma_0 |$
vary when we drop the hypothesis that $\Sigma_0$ is a set of very general points on $S$.
To get some effective results, we can apply the techniques in \cite{BS}, for smooth surfaces,
and their generalizations in \cite{KM}, for normal surfaces. These techniques are both
generalizations of Reider's theorem (see \cite{R}).

\noindent
We first recall some standard and useful definitions (see (0.1) in \cite{BS}).

\begin{definition}\label{def:23}
Let $S$ be a smooth projective surface. A line bundle $L$ on $S$ is said to be
$k$-{\em very ample}, $k \geq 0$, if the restriction map $$H^0(S,L) \to H^0(Z, {\Oc}_Z(L))$$is
surjective for any $Z \in {\mathcal Hilb}^{k+1}(S)$.
\end{definition}As an immediate consequence of the definition, we have the following
generalization of Proposition \ref{prop:9} to the surface case.

\begin{lemma}\label{lem:24}
If $L_1, \cdots, \; L_k$  are very ample line bundles on $S$, $L_1 \otimes \cdots \otimes
L_k$ is $k$-very ample.
\end{lemma}
\begin{proof}
See Lemma 0.1.1 in \cite{BS}.
\end{proof}

\begin{coroll}\label{cor:25}
Let $X \subset \PR$ be a smooth threefold and let ${\Oc}_X(1)$ denote its hyperplane bundle.
Let $\E$ be a globally generated rank-two vector bundle on $X$ and let
$[s] \in {\V}_{\delta}({\E}(k))$, where $k$ and $\delta$ are positive integers. Let $C = V(s)$,
$\Sigma = Sing(C)$ and assume that $\Sigma_0 \subseteq \Sigma$ lies on a smooth divisor
$S \subset X$, whose hyperplane
section is denoted by $H_S$. Assume further that
$$h^1(X, {\Oc}_X(k))= h^1(X, {\Ii}_{S/X}(k)) =  h^2(X, {\Ii}_{S/X}(k)) = 0.$$Then,
if $\delta_0 \leq k+1$, we have $h^1(X, {\Ii}_{\Sigma_0/X}(k)) = 0$. Therefore, if
$[s] \in {\V}_{\delta}({\E}(k))$ fails to be a regular point, the failure of the regularity
property depends on the nodes in $\Sigma \setminus \Sigma_0$. In particular, if
$\Sigma = \Sigma_0$, then $[s] \in {\V}_{\delta}({\E}(k))$ is regular.
\end{coroll}
\begin{proof}
Since ${\Oc}_S(kH_S) = {\Oc}_S(H_S)^{\otimes k}$ and since ${\Oc}_S(H_S)$ is very ample on $S$,
then ${\Oc}_S(kH_S)$ is $(\delta_0 -1)$-very ample if $k \geq \delta_0 -1$. One concludes
by using Proposition \ref{prop:16} and Remark \ref{rem:17}.
\end{proof}

\begin{remark}\label{rem:26}
\normalfont{
Observe that, since there is no assumption on the postulation of
points in $\Sigma_0$, we refind
a uniform and linear upper-bound on $\delta_0$ as we determined in Proposition \ref{prop:9} and
in Theorem \ref{thm:9bis} for $\Sigma \subset X$, where $X$ a smooth, projective threefold.
}
\end{remark}

\begin{proposition}\label{prop:27}
Let $X\subset \PR$ be a smooth threefold of degree $d$ and let $\Oc_X(1)$ denote
its hyperplane bundle. Assume that $\omega_X \cong {\Oc}_X(-m)$ for some non-negative
integer $m$. Let $\E$ be a globally generated rank-two vector bundle on $X$,
$k$ and $\delta$ be positive integers and $[s] \in {\V}_{\delta}({\E}(k))$. Let $C = V(s)$,
$\Sigma = Sing (C)$ and assume that $\Sigma_0 \subseteq \Sigma$ lies on a smooth irreducible
divisor $S \sim aH$ on $X$. Denote by $H_S$ the hyperplane section of $S$. Assume further
that$$(*) \;\;\; h^1(X, {\Oc}_X(k-a)) = h^1(X, {\Oc}_X(k))= h^2(X, {\Oc}_X(k-a)) =0.$$If $k +m > a $
and if $\delta_0 \leq \frac{k+m-a}{2} deg(D)$, for each curve $D$ on $S$,
then $h^1(X, {\Ii}_{\Sigma_0/X}(k) )=0$.

\noindent
If, in particular, $NS(S) \cong \ZZ[H_S]$, then $h^1(X, {\Ii}_{\Sigma_0/X}(k) )=0$ when
$\delta_0 \leq \frac{ad(k+m-a)}{2}$.

\noindent
Therefore, if $[s] \in {\V}_{\delta}({\E}(k))$ fails to be a regular point, the failure of the
regularity property depends on the behaviour of the nodes in $\Sigma \setminus \Sigma_0$.
In particular, if $\Sigma = \Sigma_0$, then $[s] \in {\V}_{\delta}({\E}(k))$ is a regular point.
\end{proposition}
\begin{proof}
Consider $L = {\Oc}_S( kH_S - K_S) = {\Oc}_S(k+m-a)$. Since $k+m >a$, $L$ is very ample
on $S$ and $LD = (k+m-a) deg(D)$, for each curve $D$ on $S$.
By Theorem 2.1 and Corollary 2.3 in \cite{BS}, if $\delta_0 \leq \frac{k+m-a}{2} deg(D)$,
then $$H^0(S, {\Oc}_S(k)) = H^0(S, {\Oc}_S(K_S +L)) \to \!\! \to H^0 ({\Oc}_{\Sigma_0}),$$for
$\Sigma_0 \in {\mathcal Hilb}^{\delta_0}(S)$. Since $L$ is very ample on $S$, by the Kodaira
vanishing and by the surjectivity above,
we have that $h^1(S, {\Ii}_{\Sigma_0/S}(k)) = 0$. From our assumption $(*)$, it follows that
$h^1(X, {\Ii}_{\Sigma_0/X}(k)) = 0$.
\end{proof}

\begin{remark}
\normalfont{
Observe that Proposition \ref{prop:27} applies, in particular, to
general surfaces in $\Pt$ of degree $a \geq 4$.
}
\end{remark}

By adapting the procedure of Theorem 3 in \cite{KM} to our situation, one can easily
extends to the case of $S$ a normal surface and prove analogous results. For brevity sake,
the interested reader is referred to \cite{KM}.

We conclude this section by studying the case when, given $[s] \in {\V}_{\delta}({\E}(k))$ and
$C= V(s)$, some nodes of $C$ are assumed to be on a given curve
$\Gamma \subset X$. As expected, we find some linear bounds in $k$ for the number of
admissible nodes of $C$ lying on $\Gamma$ in order that ${\Oc}_X(k)$ separates such points
on $\Gamma$. Precisely, we have:

\begin{theorem}\label{thm:34}
Let $X \subset \PR$ be a smooth threefold and 
let ${\Oc}_X(1)$ denote its hyperplane bundle. Let $\E$ be a globally generated rank-two 
vector bundle on $X$, $k$ and $\delta$ be positive integers. 
Consider $[s] \in {\V}_{\delta}({\E}(k))$ and let $C= V(s)$. Take 
$\Sigma_0 \subseteq \Sigma = Sing(C)$, $|\Sigma_0| = \delta_0$, a (not necessarily proper) 
subset of its 
nodes. Assume that $\Sigma_0$ 
lies on a local complete intersection curve $\Gamma \subset X$ such that its dualizing 
sheaf is $\omega_{\Gamma} \cong {\Oc}_{\Gamma}(e)$, for some $e \in \ZZ$. 
Assume that $h^1(\Oc_X(k)) = 0$ and that:

\noindent
(i) $k > e$,

\noindent
(ii) $h^1(X, {\Ii}_{\Gamma/X}(k))=0$, and

\noindent
(iii) $\delta_0 < deg({\Oc}_{\Gamma}(k-e))$. 

\noindent
Then $h^1(X, {\Ii}_{\Sigma_0/X}(k)) = 0$. Therefore, if 
$[s] \in {\V}_{\delta}({\E}(k))$ fails to be a regular point, the failure depends on the 
behaviour of the nodes in $\Sigma \setminus \Sigma_0$. In particular, if $\Sigma = \Sigma_0$, 
then $[s] \in {\V}_{\delta}({\E}(k))$ is a regular point.  
\end{theorem}
\begin{proof}
Consider the ideal sequence$$0 \to {\Ii}_{\Gamma/X}(k) \to {\Ii}_{\Sigma_0/X}(k) 
\to {\Ii}_{\Sigma_0/\Gamma}(k) \to 0.$$Since $h^1(X, {\Ii}_{\Gamma/X}(k))=0$, then 
$$h^1({\Ii}_{\Sigma_0/X}(k)) \hookrightarrow h^1({\Ii}_{\Sigma_0/\Gamma}(k)).$$A 
sufficient condition will be therefore the vanishing $h^1({\Ii}_{\Sigma_0/\Gamma}(k)) = 0$. 
As in \cite{G}, consider that $H^1({\Ii}_{\Sigma_0/\Gamma}(k))^{\vee} \cong 
{\rm Hom} ({\Ii}_{\Sigma_0/\Gamma}(k)), \omega_{\Gamma})$. 

Assume, by contradiction, that $h^1(X, {\Ii}_{\Sigma_0/X}(k)) \neq 0$; then a non 
zero-element corresponds to a non-zero sheaf morphism 
$${\Ii}_{\Sigma_0/\Gamma}(k)\stackrel{\varphi}{\to} 
\omega_{\Gamma}.$$Since ${\Ii}_{\Sigma_0/\Gamma}(k) \subset {\Oc}_{\Gamma}(k)$ is a 
torsion-free (but not locally free) sheaf on $\Gamma$ of rank one, $\varphi$ is an injective 
sheaf morphism. Thus, 
$$0 \to {\Ii}_{\Sigma_0/\Gamma}(k) \to \omega_{\Gamma} \to \underline{t} \to 0,$$where 
$\underline{t}$ is a torsion sheaf on $\Gamma$. Therefore, if $\chi(-)$ denotes the Euler 
characteristic,  
$$0 \leq \chi(\underline{t}) = \chi(\omega_{\Gamma}) - \chi({\Ii}_{\Sigma_0/\Gamma}(k)) = 
p_a(\Gamma) - 1 - \chi({\Oc}_{\Gamma}(k)) + \delta_0,$$i.e. 
$\delta_0 \geq \chi({\Oc}_{\Gamma}(k)) -p_a(\Gamma) + 1 = deg({\Oc}_{\Gamma}(k -e)) > 0$, since 
$k >e$ by assumption. Therefore, if $h^1({\Ii}_{\Sigma_0/\Gamma}(k)) \neq 0$, 
$\delta_0$ must be greater than or equal to $deg({\Oc}_{\Gamma}(k -e))$, which contradicts our 
hypotheses.
\end{proof}

\begin{coroll}\label{cor:35}
Let $\E$ be a globally generated rank-two 
vector bundle on $\Pt$, $k$ and $\delta$ be positive integers. 
Consider $[s] \in {\V}_{\delta}({\E}(k))$ and let $C=V(s)$. Take 
$\Sigma_0 \subseteq \Sigma = Sing(C)$ a subset of the 
nodes of $C$ s.t $|\Sigma_0| = \delta_0$. Assume that $\Sigma_0$ 
lies on a complete intersection curve $\Gamma_b \subset \Pt$ of degree $b$ 
such that $\omega_{\Gamma_b} \cong {\Oc}_{\Gamma_b}(e)$, for some $e \in \ZZ$.
Assume that $k > e$ and that $\delta_0 < b(k-e)$. Then $h^1(X, {\Ii}_{\Sigma_0/{\Pt}}(k)) = 0$. 

\noindent
Therefore, if 
$[s] \in {\V}_{\delta}({\E}(k))$ fails to be a regular point, the failure depends on the 
behaviour of the nodes in $\Sigma \setminus \Sigma_0$. In particular, if $\Sigma = \Sigma_0$, 
then $[s] \in {\V}_{\delta}({\E}(k))$ is a regular point.  

\end{coroll}
\begin{proof}
Since $\Gamma_b$ is a complete intersection, then $h^1({\Ii}_{\Gamma/\Pt}(k))=0$, for each $k$. 
Moreover, $deg({\Oc}_{\Gamma_b}(k-e)) = b(k-e)$. 
\end{proof}

\begin{remark}\label{rem:36}
\normalfont{
Observe that, with Theorem \ref{thm:34} and Corollary \ref{cor:35}, one can very easily 
construct several examples of curves on a threefold $X$ or, more specifically, 
in $\Pt$ which correspond to non-regular points of ${\V}_{\delta}({\E}(k))$, 
proving the almost-sharpness of our bounds. For example, Corollary \ref{cor:35} gives a complete 
generalization of Example 3.2 in \cite{BC}. The authors considered ${\E} = {\Oc}_{\Pt}(1) 
\oplus {\Oc}_{\Pt}(4)$ and $\delta = k+4$, with $k >>0$; they constructed 
a curve $C \subset \Pt$, corresponding to a point $[s] \in {\V}_{k+4}({\E}(k))$, with
$\Sigma_0 = \Sigma = Sing(C)$ is the set of its $(k+4)$ 
nodes lying on a line $L$ and they showed that
$[s] $ is not a regular point of ${\V}_{\delta}({\E}(k))$. Such an example is a particular  
case of our result; indeed, to hope the regularity of $[s]$ one should 
impose that the number of nodes lying on $L$ must be
$\delta_0 < 1(k-(-2)) = k+2$. Observe also that our Corollary \ref{cor:35} holds
not only for $k>>0$ but for effective values of $k$ and, furthermore, that 
we can substitute the line 
$L$ with any other complete intesection curve in $\Pt$.
}
\end{remark}

\section{Some geometric properties of space curves parametrized by
${\V}_{\delta}({\E}(k))$}\label{S:6}

From now on, we shall consider only curves in $\Pt$. Let $[s_0] \in {\V}_{\delta}({\E}(k))$,
where $\E$ is a globally generated rank-two vector bundle on $\Pt$ and where $k$ and
$\delta$ are positive integers. Let $C=V(s_0) \subset \Pt$ and let
$\Sigma$ denote the set of nodes of $C$.

The aim of this section is to study some interesting geometric properties of
the curves determined by elements in ${\V}_{\delta}({\E}(k))$. Precisely,
given
an integer $h > 0$, we can consider
smooth curves $\Gamma_h \subset \Pt$, which are defined as zero-loci of
suitable global sections $s_h \in H^0({\Pt}, {\E}(k+h))$; then, we want to study
geometric properties of the pairs $(C, \Gamma_h)$, $h \in \NN$, expecially from the liaison relation
point of view. From Rao's paper (\cite{Rao}), we know that $C$ and $\Gamma_h$ lie in the
same liaison class, for each $h >0$. Indeed, since such curves correspond to global sections of
twists of the same vector bundle $\E$, by the Koszul exact sequences in $\Pt$, their
Rao's modules are isomorphic up to the shift of $h$, i.e.
\begin{equation}\label{eq:rao}
M(C) = \bigoplus_t H^1(C, {\Ii}_{C/{\Pt}}(t)) \cong \bigoplus_t H^1(\Gamma_h,
{\Ii}_{C/{\Pt}}(h + t)) = M(\Gamma_h).
\end{equation}One can be more precise by recalling
the following terminology.

\begin{definition}\label{def:38}(see Def. 1.1 in \cite{Rao})
Let $V_1, \; V_2 \subset \Pt$ be two subschemes which are locally Cohen-Macaulay and
equidimensional of codimensione two. Denote by $\sim_l$ the relation of
(geometric) linkage in $\Pt$. Then
$V_1$ and $V_2$ are in the same {\em even liaison class} if
$$V_1 \sim_l Z_1 \sim_l \cdots \sim_l Z_{2k+1} \sim_l V_2,$$for some schemes $Z_i \subset \Pt$.
In such a case, $V_1$ and $V_2$
are said to be {\em evenly linked}. The resulting equivalence relation is called
{\em biliaison}.
\end{definition}By (\ref{eq:rao}) and by Lemma 1.6 in \cite{Rao}, we immediately
observe that the curves $C$ and $\Gamma_h$ lie in the
same biliaison class, apart from some particular cases. Indeed, since $C$ and $\Gamma_h$ are
both subcanonical curves in $\Pt$, it may happen that $C$ and $\Gamma_h$ are directly (so oddly)
linked in $\Pt$ (see \cite{C}). However, when this last situation occurs, there are strong
restrictions not only on the vector bundle $\E$ but also on the numerical characters of
$\Gamma_h$ (degree, genus and postulation); precisely, one has at most three possible cases
for $\Gamma_h$ (for details, the reader is referred to the original paper \cite{C}).

Here we are interested in analyzing the geometric properties that $C$ and
$\Gamma_h$ share, for each $h >0$, particularly towards the biliaison relation between $C$ and
$\Gamma_h$ (i.e. with general choices of $\E$ and $\Gamma_h$).

\noindent
By recalling Remark \ref{rem:pic}, a first result is the following.

\begin{theorem}\label{thm:37}
Let $\E$ be a globally generated rank-two vector bundle on $\Pt$ and denote by
$c_1$ the first Chern class of $\E$. Let $k$ and
$\delta$ be positive integers such that $\delta \leq k+1$. Let $[s_0] \in {\V}_{\delta}({\E}(k))$ 
correspond to an irreducible, nodal curve $C\subset \Pt$, whose set of nodes is 
denoted by $\Sigma$. Take $h \geq 1$ an integer such that
\begin{equation}\label{eq:k+h}
k+h \geq  c_1 \in \ZZ.
\end{equation}
Then, there always exist:
\begin{itemize}
\item[(i)] a smooth, irreducible curve $\Gamma_h \subset \Pt$ simply passing through
$\Sigma$ and such that $\Gamma_h = V(s_h)$, where $[s_h] \in {\Pp}(H^0( {\E}(k+h))$,
\item[(ii)] a normal surface $S_{k,h} \subset \Pt$ of degree $d_{k,h} = c_1 + 2k + h$
containing both $C$ and $\Gamma_h$.
\end{itemize}
\end{theorem}
\begin{proof}
As in Theorem \ref{prop:3.fundamental}, we consider the smooth projective fourfold
${\mathcal P} = {\Pp}_{\Pt}({\E}(k))$,
together with the surjective morphism $\pi :  {\mathcal P} \to \Pt$ and with the tautological
line bundle ${\Oc}_{\mathcal P}(1)$. Thus the curve $C$, which is the zero-locus of the
global section $s_0 \in H^0({\E}(k))$, corresponds to an effective
divisor $D_0 \in |{\Oc}_{\mathcal P}(1) |$.

\noindent
(i) Given $h >0$, we first want to construct a section
$s_h \in H^0({\E}(k+h))$, which corresponds to a divisor
$D_h \in | {\Oc}_{\mathcal P}(1) \otimes \pi^* ({\Oc}_{\Pt}(h))|$ simply vanishing along the
lines of the scheme ${\La} = \bigcup_{i=1}^{\delta} L_i := \pi^{-1}(\Sigma)$. Such
a divisor will correspond to the curve $\Gamma_h$  we want to determine.

Recall that the divisor $D_0$ corresponding to $C$
is singular having $\delta$ rational double points
at $\Sigma^1 \subset {\La}$, where $\Sigma^1 \cong \Sigma$
(see the proof of Theorem \ref{prop:3.fundamental}).
Assume for a moment that the sheaf
\begin{equation}\label{eq:sheaf}
{\Ii}_{{\La}/{\Pmc}} \otimes {\Oc}_{\mathcal P}(1) \otimes \pi^* ({\Oc}_{\Pt}(h))
\end{equation}is globally generated on $\Pmc$
(we shall show this fact later on in this proof); thus, the scheme $\La$ coincides with the base locus
of the linear system
$|{\Ii}_{{\La}/{\Pmc}} \otimes {\Oc}_{\mathcal P}(1) \otimes \pi^* ({\Oc}_{\Pt}(h))|$. Since
$H^0({\E}(k)) \hookrightarrow H^0({\E}(k+h))$, for a general
$\sigma \in H^0({\Oc}_{\Pt}(h))$ the global section $s_0 \otimes \sigma \in H^0({\E}(k+h))$
behaves as $s_0$ around the points in $\Sigma$; moreover, since
$\pi_*({\Oc}_{\mathcal P}(1) \otimes \pi^* ({\Oc}_{\Pt}(h))) \cong {\E}(k+h)$, if
$$B := {\rm Base \; scheme}(|{\Ii}_{{\La}/{\Pmc}}\otimes {\Oc}_{\mathcal P}(1)
\otimes \pi^* ({\Oc}_{\Pt}(h))|)$$then ${\La}=B$ as schemes, which means that each
line $L_i$ is reduced in $B$.

We want to show that the general element of
$|{\Ii}_{{\La}/{\Pmc}} \otimes {\Oc}_{\mathcal P}(1) \otimes \pi^* ({\Oc}_{\Pt}(h))|$
is smooth along ${\La}$. Denote by $D_h$ the general divisor of
$|{\Ii}_{{\La}/{\Pmc}} \otimes {\Oc}_{\mathcal P}(1) \otimes \pi^* ({\Oc}_{\Pt}(h))|$;
since we assumed that the sheaf in (\ref{eq:sheaf}) is globally generated, then
$${\N}_{{\La}/D_h}^{\vee} \otimes {\Oc}_{\mathcal P}(D_h)
\cong \frac{{\Ii}_{{\La}/{\Pmc}}}{{\Ii}_{{\La}/{\mathcal P}}^2} \otimes
{\Oc}_{\mathcal P}(1) \otimes
\pi^* ({\Oc}_{\Pt}(h))$$is globally generated. Since
we have ${\La} \subset D_h \subset {\mathcal P}$, we get
\begin{equation}\label{eq:37**}
0 \to {\Oc}_{D_h} \cong {\N}_{D_h/{\mathcal P}}^{\vee} \otimes {\Oc}_{\mathcal P}(D_h)
\to {\N}_{{\La}/{\mathcal P}}^{\vee} \otimes {\Oc}_{\mathcal P}(D_h)  \to
{\N}_{{\La}/D_h}^{\vee} \otimes {\Oc}_{\mathcal P}(D_h) \to 0.
\end{equation}Fix $L= L_{i_0} \subset {\La}$, for some $1 \leq i_0 \leq \delta$,
and restrict (\ref{eq:37**}) to $L$; therefore we have
\begin{equation}\label{eq:37***}
{\Oc}_{L} \to {\N}_{{\La}/{\mathcal P}}^{\vee} \otimes {\Oc}_{L}(D_h)  \to
{\N}_{{\La}/D_h}^{\vee} \otimes {\Oc}_{L}(D_h) \to 0.
\end{equation}Since ${\N}_{{\La}/{\mathcal P}}^{\vee} \otimes {\Oc}_{L}(D_h)$
is a globally generated rank-three vector bundle on the line $L$, there exists a global
section nowhere vanishing on $L$; this implies that (\ref{eq:37***}) is exact
and that ${\N}_{{\La}/{\mathcal P}}^{\vee} \otimes {\Oc}_{L}(D_h)$ is locally
free on $L$. So it is ${\N}_{{\La}/{\mathcal P}}^{\vee}|_L$. Since we have
$L \subset {\La} \subset D_h$, then also ${\N}_{L/D_h}^{\vee}$ is locally free. Therefore,
since $L$ is smooth in ${\mathcal P}$, we have

\begin{displaymath}
\begin{array}{ccccccc}
 & & & 0 & & & \\
 & & & \downarrow & & & \\
  & {\N}_{D_h/{\mathcal P}}^{\vee}|_L & \to & \Omega^1_{{\mathcal P}|_{D|_L}} &
\to & \Omega^1_{D_h}|_L & \to 0 \\
 & \downarrow & & \shortparallel & & \downarrow^{\alpha} & \\
0 \to & {\N}_{L/{\mathcal P}}^{\vee} & \to & \Omega^1_{{\mathcal P}|_L} & \to &
\Omega^1_L \cong
{\Oc}_{{\Pp}^1}(-2) & \to 0 \\
 & \downarrow & & \downarrow & &  \downarrow^{\;} & \\
 & {\N}_{L/D_h}^{\vee}& & 0 & &  0 & \\
 & \downarrow & &  & & &\\
 & 0 & &  & & &
 \end{array}
\end{displaymath}By the Snake lemma, $ker(\alpha) \cong {\N}_{L/D_h}^{\vee}$ is locally free
on $L$. This implies that $\Omega^1_{D_h}|_L $ is locally free on $L$, i.e. the general
element of $|{\Ii}_{{\La}/{\Pmc}} \otimes {\Oc}_{\mathcal P}(1)
\otimes \pi^* ({\Oc}_{\Pt}(h))|$ is smooth on $L$.

To complete the proof, we only have to show that the sheaf in (\ref{eq:sheaf}) is actually
globally generated on $\Pmc$.  Observe that
${\Ii}_{{\La}/{\Pmc}} \otimes {\Oc}_{\mathcal P}(1) \otimes \pi^* ({\Oc}_{\Pt}(h))$ is
globally generated iff ${\Ii}_{\Sigma/\Pt} \otimes {\E}(k+h)$ is globally generated
on $\Pt$. Since a sufficient condition for the global generation
of ${\Ii}_{\Sigma/\Pt} \otimes {\E}(k+h)$ is its $0$-regularity as a sheaf on $\Pt$,
observe that ${\Ii}_{\Sigma/\Pt} \otimes {\E}(k+h)$ is $0$-regular iff
\begin{equation}\label{eq:37****}
h^1({\Ii}_{\Sigma/\Pt} \otimes {\E}(k+h-1)) = h^2({\Ii}_{\Sigma/\Pt} \otimes {\E}(k+h-2))
= h^3({\Ii}_{\Sigma/\Pt} \otimes {\E}(k+h-3))=0.
\end{equation}From the Griffiths vanishing result (see \cite{SS}, page 107), it follows that
\begin{equation}\label{eq:37*****}
h^1({\E}(k+h-1)) = h^2( {\E}(k+h-2)) = h^3({\E}(k+h-3))=0,
\end{equation}since $\E$ is globally generated and since $k+h \geq c_1 = det({\E})$
by assumption. Therefore, from the exact sequence
$$0 \to {\Ii}_{\Sigma/\Pt} \otimes {\E}(k+h) \to {\E}(k+h) \to {\E}(k+h)|_{\Sigma} \to 0$$and from
(\ref{eq:37*****}), the last two equalities in
(\ref{eq:37****}) hold. It only remains to
show that $h^1({\Ii}_{\Sigma/\Pt} \otimes {\E}(k+h-1)) =0$. To this aim, take
\[
\begin{array}{llllcl}
0 \to & H^0({\Ii}_{\Sigma/\Pt} \otimes {\E}(k+h-1)) & \to & H^0({\E}(k+h-1)) &
\stackrel{\alpha_{k+h-1}}{\longrightarrow} & H^0({\E}(k+h-1)|_{\Sigma})\\
 & & \to & H^1({\Ii}_{\Sigma/\Pt} \otimes {\E}(k+h-1))& \to & 0.
\end{array}
\]Since $\delta \leq k+1$ by assumption, then
$$H^0({\E}(k)) \stackrel{\alpha_{k}}{\longrightarrow} H^0({\E}(k)|_{\Sigma})$$is surjective;
therefore, $\alpha_{k+h-1}$ is surjective since $h \geq 1$, i.e.
$ H^1({\Ii}_{\Sigma/\Pt} \otimes {\E}(k+h-1)) = (0)$.

\noindent
(ii) After having constructed the curve $\Gamma_h$, which corresponds to the general element
of $|{\Ii}_{{\La}/{\Pmc}} \otimes {\Oc}_{\mathcal P}(1) \otimes \pi^* ({\Oc}_{\Pt}(h))|$, take
$s_h \in H^0({\E}(k+h))$ such that $\Gamma_h = V(s_h)$. We can consider the rank-two vector
bundle morphism
$$\tau = (s_0 , s_h) : {\Oc}_{\Pt} \oplus {\Oc}_{\Pt}(-h) \to {\E}(k),$$where
$s_0 \in H^0({\E}(k))$ is such that $C = V(s_0)$. The degeneration locus of the morphism $\tau$
is a surface $S_{k,h} = V(det(\tau))$, where $det(\tau) \in
H^0({\Pt}, \; det({\E}) \otimes {\Oc}_{\Pt} (2k+h))$. Therefore, $S_{k,h}$ is a
surface of degree $d_{k,h} = c_1 + 2k+h$ containing both $C=V(s_0)$ and
$\Gamma_h=V(s_h)$. Its singular locus is determined by the condition $rank(\tau) < 1$; therefore,
by the construction of $\Gamma_h$ in (i), we immediately observe that $\Sigma \subseteq
Supp(C \cap \Gamma_h) = Sing(S_{k,h})$.
\end{proof}For each $k$ and $h$ as in \eqref{eq:k+h},
the surface $S_{k,h}$ constructed in Theorem \ref{thm:37} is the most "natural" surface in
$\Pt$ containing both curves $C$ and $\Gamma_h$. One can very easily deduce some biliaison properties of
$C$ and $\Gamma_h$ on $S_{k,h}$. We first recall the following more general
definition from \cite{H}.

\begin{definition}\label{def:40}
Let $V_1$, $V_2$ be schemes of equidimension one without embedded components
which are evenly linked in $\Pt$. $V_2$ is said to be
obtained from $V_1$ by an {\em elementary biliaison of height} $h$, for some $h \in \ZZ$, if there exists
a surface
$S \subset \Pt$ containing $V_1$ and $V_2$ such that $V_2 \sim V_1 + hH$ as generalized divisor
on $S$, where $H$ is the plane section of $S$. This is equivalent to saying there exist
surfaces $T_1$, of degree $t_1$ containing $V_1$, and $T_2$, of degree $t_2$ containing $V_2$,
such that $t_2 = t_1 +h $ and the scheme $W_1$ linked to $V_1$ by $T_1 \cap S$ is equal to the scheme
$W_2$ linked to $V_2$ by $T_2 \cap S$ (see Proposition 4.3 (b) in \cite{H} or \cite{LR}
pg. 276)
\end{definition}Thus, we can state the following:

\begin{proposition}\label{prop:41}
Let $[s_0] \in {\V}_{\delta}({\E}(k))$, with $\E$ a non-splitting, globally generated
rank-two vector bundle on $\Pt$ and with $k$ and $\delta$ positive integers such that $\delta \leq k+1$. 
Let $C= V(s_0)$ and $\Sigma = Sing(C)$.
Then, for each
integer $h\geq 1$ as in \eqref{eq:k+h}, each curve $\Gamma_h$ as in Theorem \ref{thm:37} (i)
is obtained from $C$ by an elementary biliaison of height $h$ on the surface $S_{k,h}$
\end{proposition}
\begin{proof}For simplicity of notation, we denote by $S$ the surface $S_{k,h}$ and by 
$c_1$ the first Chern class of $\E$. Denote by $\G$ the cokernel of the map
$$\tau = (s_0 , s_h) : {\Oc}_{\Pt} \oplus {\Oc}_{\Pt}(-h) \to {\E}(k),$$where
$s_0 \in H^0({\E}(k))$ and $s_h \in H^0({\E}(k+h))$ are
such that $C = V(s_0)$ and $\Gamma_h = V(s_h)$. By the diagram
\begin{equation}\label{eq:(*6)}
\begin{aligned}
\begin{array}{rcccclr}
 & & &  & & 0 & \\
 & & &  & & \downarrow & \\
0 \to  & {\Oc}_{\Pt}  &\to & {\Oc}_{\Pt} \oplus {\Oc}_{\Pt}(-h) &\to & {\Oc}_{\Pt}(-h)  & \to 0 \\
 & \shortparallel & & \downarrow^{\tau} & & \downarrow^{det(s_h)} &  \\
0 \to  & {\Oc}_{\Pt}  & \stackrel{s_0}{\to} & {\E}(k) &\to & {\Ii}_{C/\Pt}(c_1+2k) & \to 0 \\
 & & & \downarrow & & \downarrow & \\
   &  &  & {\G} &   & {\Ii}_{C/S}(c_1+2k) &  \\
 & & & \downarrow & & \downarrow & \\
 &  & & 0 & & 0 &,
\end{array}
\end{aligned}
\end{equation}we see that $\G \cong {\Ii}_{C/S}(c_1)$. Reversing the roles of $s_0$ and $s_h$ in
diagram \eqref{eq:(*6)}, we similarly find that $\G \cong {\Ii}_{\Gamma_h /S}(c_1+h)$. Hence, we get
$ {\Ii}_{C/S} \cong  {\Ii}_{\Gamma_h /S}(h)$; since $\E$ 
is non-splitting, neither $C$ nor $\Gamma_h$ can be equivalent to multiples of the hyperplane sections 
of $S$. Thus, by Definition \ref{def:40},
$\Gamma_h$ is obtained by an elementary biliason of height $h$ on $S$.
\end{proof}

\end{document}